\documentclass[11pt,draft]{article}
\usepackage[letterpaper]{geometry}

\usepackage{amsmath,amsthm,amssymb,mathrsfs}
\usepackage{colonequals}
\usepackage{enumerate}
\usepackage{graphics}
\usepackage[all]{xy}
\usepackage{tikz} 
\usetikzlibrary{arrows,decorations.pathmorphing,backgrounds,positioning,fit}
\usepackage{fixme}
\usepackage{citesort}
\usepackage{graphicx}
\DeclareGraphicsExtensions{.pdf,.png,.jpg}

\newcommand{\eqclass}[1]{[\![#1]\!]}

\newcommand{\modelt}{\models^+}
\newcommand{\modelf}{\models^-}

\newcommand{\concat}{{^{\frown}}}

\DeclareMathOperator{\subf}{Subf}
\DeclareMathOperator{\lit}{Lit}
\DeclareMathOperator{\Heads}{Heads}
\DeclareMathOperator{\Tails}{Tails}
\DeclareMathOperator{\Monday}{Monday}
\DeclareMathOperator{\Tuesday}{Tuesday}
\DeclareMathOperator{\Awake}{Awake}
\newcommand{\hor}[1]{\lor\!\!_{/#1}\,}
\newcommand{\hand}[1]{\land\!_{/#1}\,}

\newcommand{\abelard}{Abelard}
\newcommand{\eloise}{Eloise}
\let \str = \mathbf
\DeclareMathOperator{\dom}{dom}

\newcommand\stochastic[0]{\reflectbox{\ensuremath{\mathsf{S}}}} 

\theoremstyle{definition}
\newtheorem*{defn}{Definition}

\theoremstyle{remark}

\newtheorem*{example}{Example}

\let	\phi	=	\varphi

\newcommand{\vaananen}{V\"a\"an\"anen}

\title{A logical analysis of Monty Hall and Sleeping Beauty}
\author{Allen L. Mann\footnote{The first author wishes to gratefully acknowledge the partial support of the European Science Foundation EUROCORES program LogICCC [FP002---Logic for Interaction (LINT)] and the Academy of Finland (grant 129208).}\ \ and Ville Aarnio}

\begin{document}
\maketitle


\begin{abstract}
Hintikka and Sandu's independence-friendly (IF) logic is a conservative extension of first-order logic that allows one to consider semantic games with imperfect information. 
In the present article, we first show how several variants of the Monty Hall problem can be modeled as semantic games for IF sentences. 
In the process, we extend IF logic to include semantic games with chance moves and dub this extension \emph{stochastic IF logic}. 
Finally, we use stochastic IF logic to analyze the Sleeping Beauty problem, leading to the conclusion that the thirders are correct while identifying the main error in the halfers' argument.
\end{abstract}


\section{Introduction}
A game-show contestant is presented with three doors, 
one of which conceals a prize. 
After the contestant selects a door,
the host (Monty Hall) opens one of the two remaining doors,
being careful not to reveal the prize.
He then offers the contestant 
the opportunity to switch doors. 
Should the contestant stick with her original choice or switch?%
\footnote{The Monty Hall problem was popularized by Marilyn vos Savant \cite{Savant:1990}. 
For an overview of the history of the problem and its wider implications, see Tierney \cite{Tierney:1991,Tierney:2008}.}

The surprising answer is that the contestant should always switch doors.
If she does, 
she will win with probability 2/3,
whereas 
if she sticks with her original door 
she will win with probability 1/3.
What is even more surprising is that 
the Monty Hall problem 
can be expressed by the independence-friendly sentence
\(\varphi_{\mathrm{MH}}\),
\[
	\forall x 
	\bigl(\exists y/\{x\}\bigr)
	\forall z
	\Bigl[\,
		(z \not= x \land z \not= y) 
		\implies 
		\bigl(\exists y/\{x\}\bigr)\,
		x = y\,
	\Bigr],
\]
interpreted in the structure 
\(\str{M} = \{1,2,3\}\).
The sentence can be read: 
\begin{quote}
	For any Door $x$ concealing the prize,  
	there exists a Door $y$ 
	chosen (independently of $x$) by the contestant 
	such that
	for every possible Door $z$ opened by Monty Hall 
	that differs from the door containing the prize 
	and
	the door chosen by the contestant,
	there is a Door $y$ chosen (independently of $x$) by the contestant 
	such that \(x = y\).
\end{quote}

To understand the connection between 
the Monty Hall problem and $\varphi_{\mathrm{MH}}$, 
one must interpret the formula game theoretically. 
The semantic game for a first-order sentence is 
a contest between two players.
The existential player attempts to verify the sentence by
choosing the values of existentially quantified variables, 
while the universal player tries to falsify the sentence 
by picking the values of universally quantified variables. 
Disjunctions prompt the verifier to choose a disjunct;
conjunctions prompt the falsifier to pick a conjunct.
Negation tells the players to switch roles. 
The sentence is true 
if the existential player has a winning strategy, 
false if the universal player has a winning strategy.
Traditionally, 
the existential player is named \eloise\
and 
universal player \abelard.

The semantic game for a first-order sentence 
is a game with perfect information
in the sense that, 
at each decision point, 
the active player is aware of all prior moves.
In contrast, 
the semantic game for 
$\varphi_{\mathrm{MH}}$ 
is a game with imperfect information
since \eloise\ must choose the value of $y$ without knowing the value of $x$.
Thus, 
the sentence 
$\varphi_{\mathrm{MH}}$ 
is ``independence-friendly'' 
in the sense that the contestant's choices may not depend on the location of the prize.

The present paper is organized as follows.
The next section 
presents the Monty Hall problem as an extensive game, 
while Section \ref{Section:IFLogic}
introduces the syntax and semantics of independence-friendly (IF) logic.
We then show that the Monty Hall problem is equivalent to the semantic game for $\varphi_{\mathrm{MH}}$ 
and consider a variant of the problem 
in which the host is not required to offer the contestant the opportunity to switch doors.
In Section \ref{Section:Stochastic IF logic},
we consider another variant of the Monty Hall problem in which the host is indifferent to the outcome of the game.
To analyze this variant, 
we extend IF logic by adding a third player (Nature) 
who makes moves at random. 
We then use this extension in Section \ref{Section:SleepingBeauty} to address the controversy surrounding the Sleeping Beauty problem,
which has divided the philosophical community into two camps, thirders and halfers, who respectively defend the contradictory solutions 1/3 and 1/2.
We briefly review the arguments presented by both sides, 
and present several ways to formalize the problem using stochastic IF logic. 
Under our preferred formalization
the answer is 1/3, 
but under two alternate formalizations 
the answer is 1/2.

\section{The Monty Hall problem as an extensive game}
\label{Section:The Monty Hall problem as an extensive game}

In an extensive game,
players take turns making moves until the game ends, 
at which point each player receives a certain payoff.
The following more formal definition is taken with slight modifications from Osborne and Rubinstein \cite[p.~200]{Osborne:1994}. 

\begin{defn} An \emph{extensive game with imperfect information} has the following components.
    \begin{itemize}
        	\item
        	A finite set $N$ of \emph{players}.
        	\item
        	A set $H$ of sequences (called \emph{histories}) closed under initial segments.
        		\begin{itemize}
        			\item 
        			A component $a_{i}$ of a history \(h = (a_{1},\dots, a_{n})\) is called an \emph{action}.
        			\item 
        			A history of the form \(h\concat a = (a_{1},\dots,a_{n},a)\) is called a \emph{successor} of $h$.
        			\item 
        			A \emph{terminal history} is a history with no successors. 
        			The set of terminal histories is denoted \(Z \subseteq H\); 
        			the set of actions available after the nonterminal history $h \in H \setminus Z$ is denoted 
        				\[
        					A(h) = \bigl\{a : h\concat a \in H\bigr\}.
        				\]
        		\end{itemize}
        	\item
        	A \emph{player function} $P\colon H\setminus Z \to N$ that indicates whose turn it is to move. 
	Let \(H_{p} = P^{-1}(p)\) denote the set of histories after which it is player $p$'s turn.
        	\item
        	For each player \(p \in N\),
	an equivalence relation $\sim_{p}$ on $H_{p}$ 
	with the property that
	for all \(h,h' \in H_{p}\),
		\[
			h \sim_{p} h' \quad\text{implies}\quad A(h) = A(h').
		\]
	If \(h \sim_{p} h'\), we say the histories $h$ and $h'$ are \emph{indistinguishable} to player $p$. 
	An equivalence class 
	\[\eqclass{h}_{\sim_{p}} = \bigl\{\,h' \in H_{p} : h \sim_{p} h'\,\bigr\}\]
	is called an \emph{information set} for player $p$.
        	\item 
        	A \emph{utility function} 
	\(u\colon Z \to \mathbb{R}^{N}\) 
	that specifies the payoff each player receives at the end of the game.
	For each terminal history \(h \in Z\) and player \(p \in N\)\!, let 
	\(u_{p}(h) = u(h)(p)\) denote the payoff received by player $p$.
	A \emph{constant-sum} game is one in which the sum of the payoffs received by the players always takes the same value.
	In a \emph{win-lose} game, 
	a single player (the \emph{winner}) receives a payoff of 1, 
	while all other players receive a payoff of 0.
    \end{itemize}
\end{defn}

\begin{figure}[htbp]
	\centering
		\resizebox{4in}{!}{
		\begin{tikzpicture}
			\node (root) 	[circle,draw] {II};
			\node (a) 		[rounded corners,draw]	
										[above right=of root, xshift=5mm, yshift=40mm] {I$^{a}$}
											edge node[auto,swap] {1} (root);
			\node (aa)		[circle,draw] 
										[above right=of a, yshift=20mm] {II}
											edge node[auto,swap] {1} (a);
			\node (aab)		[rounded corners,draw]
										[above right=of aa, yshift=-3mm]{I$^{b}$}
											edge node[auto,swap] {2} (aa);								
			\node (aac)		[rounded corners,draw]
										[below right=of aa, yshift=3mm]{I$^{c}$}
											edge node[auto] {3} (aa);
			\node (aaba)	[circle,fill] 
										[above right=of aab, yshift=-9mm, label=right:{I}]{}
											edge node[auto,swap] {1} (aab);
			\node (aabc)	[circle,fill] 
										[below right=of aab, yshift=9mm, label=right:{II}]{}
											edge node[auto] {3} (aab);
			\node (aaca)	[circle,fill] 
										[above right=of aac, yshift=-9mm, label=right:{I}]{}
											edge node[auto,swap] {1} (aac);
			\node (aacb)	[circle,fill] 
										[below right=of aac, yshift=9mm, label=right:{II}]{}
											edge node[auto] {2} (aac);
			\node (ab)		[circle,draw] 
										[right=of a, xshift=-2mm]	{II}
											edge node[auto,swap] {2} (a);
			\node (abc)		[rounded corners,draw]
										[right=of ab, xshift=-2mm]{I$^{d}$}
											edge node[auto,swap] {3} (ab);
			\node (abca)	[circle,fill] 
										[above right=of abc, yshift=-9mm, label=right:{I}]{}
											edge node[auto,swap] {1} (abc);
			\node (abcb)	[circle,fill] 
										[below right=of abc, yshift=9mm, label=right:{II}]{}
											edge node[auto] {2} (abc);
			\node (ac)		[circle,draw] 
										[below right=of a, yshift=-10mm]{II}
											edge node[auto] {3} (a);
			\node (acb)		[rounded corners,draw]
										[right=of ac, xshift=-2mm]{I$^{e}$}
											edge node[auto] {2} (ac);
			\node (acba)	[circle,fill] 
										[above right=of acb, yshift=-9mm, label=right:{I}]{}
											edge node[auto,swap] {1} (acb);
			\node (acbc)	[circle,fill] 
										[below right=of acb, yshift=9mm, label=right:{II}]{}
											edge node[auto] {3} (acb);
			\node (b)			[rounded corners,draw] 
										[right=of root, xshift=60mm]	{I$^{a}$}
											edge node[auto,swap] {2} (root);
			\node (ba)		[circle,draw] 
										[above right=of b, yshift=20mm] {II}
											edge node[auto,swap] {1} (b);
			\node (bac)		[rounded corners,draw]
										[right=of ba]{I$^{c}$}
											edge node[auto,swap] {3} (ba);
			\node (baca)	[circle,fill] 
										[above right=of bac, yshift=-9mm, label=right:{II}]{}
											edge node[auto,swap] {1} (bac);
			\node (bacb)	[circle,fill] 
										[below right=of bac, yshift=9mm, label=right:{I}]{}
											edge node[auto] {2} (bac);
			\node (bb)		[circle,draw] 
										[right=of b]	{II}
											edge node[auto,swap] {2} (b);
			\node (bba)		[rounded corners,draw]
										[above right=of bb, yshift=-3mm]{I$^{f}$}
											edge node[auto,swap] {1} (bb);								
			\node (bbc)		[rounded corners,draw]
										[below right=of bb, yshift=3mm]{I$^{d}$}
											edge node[auto] {3} (bb);
			\node (bbab)	[circle,fill] 
										[above right=of bba, yshift=-9mm, label=right:{I}]{}
											edge node[auto,swap] {2} (bba);
			\node (bbac)	[circle,fill] 
										[below right=of bba, yshift=9mm, label=right:{II}]{}
											edge node[auto] {3} (bba);
			\node (bbca)	[circle,fill] 
										[above right=of bbc, yshift=-9mm, label=right:{II}]{}
											edge node[auto,swap] {1} (bbc);
			\node (bbcb)	[circle,fill] 
										[below right=of bbc, yshift=9mm, label=right:{I}]{}
											edge node[auto] {2} (bbc);
			\node (bc)		[circle,draw] 
										[below right=of b, yshift=-20mm]{II}
											edge node[auto] {3} (b);
			\node (bca)		[rounded corners,draw]
										[right=of bc]{I$^{g}$}
											edge node[auto] {1} (bc);
			\node (bcab)	[circle,fill] 
										[above right=of bca, yshift=-9mm, label=right:{I}]{}
											edge node[auto,swap] {2} (bca);
			\node (bcac)	[circle,fill] 
										[below right=of bca, yshift=9mm, label=right:{II}]{}
											edge node[auto] {3} (bca);
			\node (c)			[rounded corners,draw] 
										[below right=of root, xshift=5mm, yshift=-40mm]{I$^{a}$}
											edge node[auto] {3} (root);
			\node (ca)		[circle,draw] 
										[above right=of c, yshift=10mm] {II}
											edge node[auto,swap] {1} (c);
			\node (cab)		[rounded corners,draw]
										[right=of ca, xshift=-2mm]{I$^{b}$}
											edge node[auto,swap] {2} (ca);
			\node (caba)	[circle,fill] 
										[above right=of cab, yshift=-9mm, label=right:{II}]{}
											edge node[auto,swap] {1} (cab);
			\node (cabc)	[circle,fill] 
										[below right=of cab, yshift=9mm, label=right:{I}]{}
											edge node[auto] {3} (cab);
			\node (cb)		[circle,draw] 
										[right=of c, xshift=-2mm]	{II}
											edge node[auto,swap] {2} (c);
			\node (cba)		[rounded corners,draw]
										[right=of cb, xshift=-2mm]{I$^{f}$}
											edge node[auto,swap] {1} (cb);
			\node (cbab)	[circle,fill] 
										[above right=of cba, yshift=-9mm, label=right:{II}]{}
											edge node[auto,swap] {2} (cba);
			\node (cbac)	[circle,fill] 
										[below right=of cba, yshift=9mm, label=right:{I}]{}
											edge node[auto] {3} (cba);
			\node (cc)		[circle,draw] 
										[below right=of c, yshift=-20mm]{II}
											edge node[auto] {3} (c);
			\node (cca)		[rounded corners,draw]
										[above right=of cc, yshift=-3mm]{I$^{g}$}
											edge node[auto,swap] {1} (cc);								
			\node (ccab)	[circle,fill] 
										[above right=of cca, yshift=-9mm, label=right:{II}]{}
											edge node[auto,swap] {2} (cca);
			\node (ccac)	[circle,fill] 
										[below right=of cca, yshift=9mm, label=right:{I}]{}
											edge node[auto] {3} (cca);
			\node (ccb)		[rounded corners,draw]
										[below right=of cc, yshift=3mm]{I$^{e}$}
											edge node[auto] {2} (cc);
			\node (ccba)	[circle,fill] 
										[above right=of ccb, yshift=-9mm, label=right:{II}]{}
											edge node[auto,swap] {1} (ccb);
			\node (ccbc)	[circle,fill] 
										[below right=of ccb, yshift=9mm, label=right:{I}]{}
											edge node[auto] {3} (ccb);
		\end{tikzpicture}
		}
			\caption{The Monty Hall problem}
  \label{Figure-MontyHallGameTree}
\end{figure}
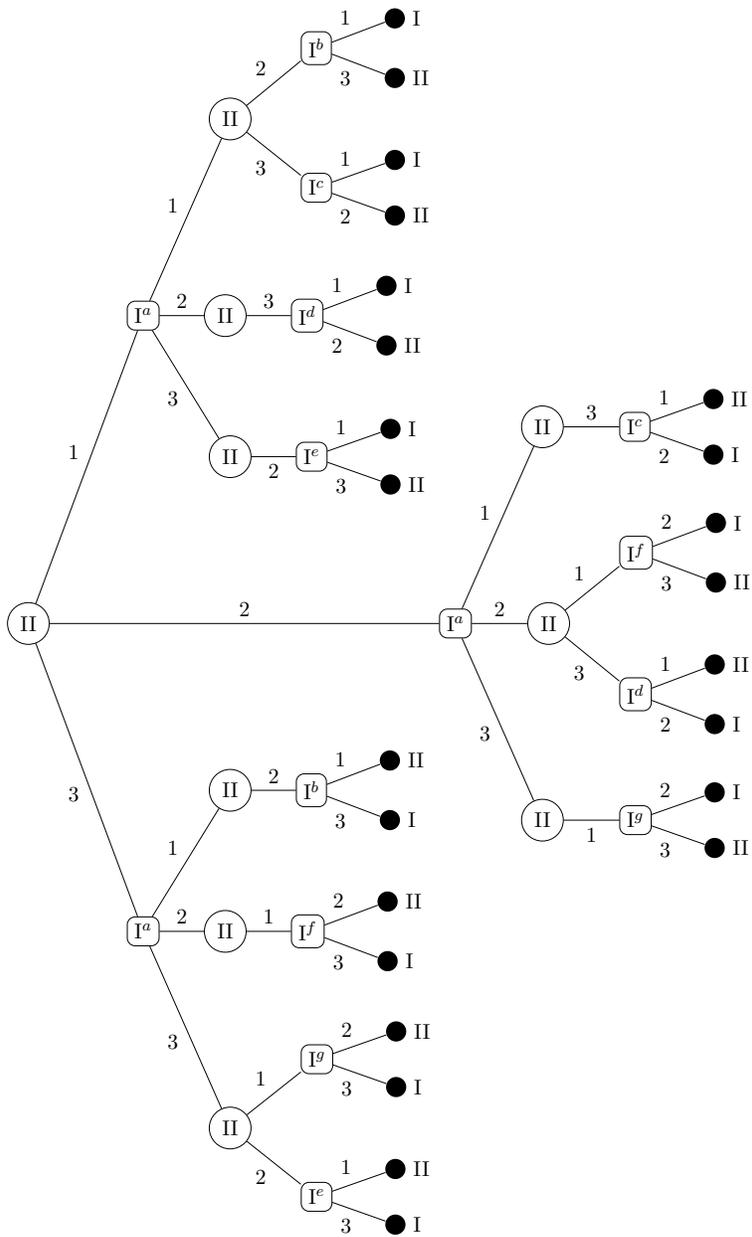

An extensive game can also be represented by a game tree \cite[page 42]{Maschler:2013}.
The game tree for the Monty Hall problem is shown in 
Figure \ref{Figure-MontyHallGameTree}, 
where for later convenience we will treat the contestant as player I, and Monty Hall as player II.
First, 
Monty Hall secretly places the grand prize behind one of three doors. 
Second, the contestant guesses which door conceals the prize.
Third, Monty opens a door, revealing its contents.
He never opens the door initially chosen by the contestant, however,
nor does he reveal the prize.
Thus, 
if the contestant guessed correctly, 
Monty Hall is free to open either of the two remaining doors. 
If the contestant guessed incorrectly, 
however,
Monty is forced to open the only other door that does not conceal the prize. 
Finally, 
Monty Hall offers the contestant the opportunity to stick with her original guess 
or to switch to the other unopened door.
The contestant wins if she chooses the door concealing the grand prize; otherwise Monty Hall wins.

What makes the Monty Hall problem interesting is the contestant's uncertainty about which door conceals the prize. 
Thus, she is unsure which node of the game tree corresponds to the current state of the game. 
In Figure \ref{Figure-MontyHallGameTree}, 
we label positions that are indistinguishable to the contestant with the same superscript letter. 
For example, 
the contestant cannot distinguish between the positions labeled I$^{a}$ because she does not know the location of the prize when making her initial choice. 
Similarly, in both positions labeled I$^{b}$\!, the contestant initially chose Door 1, after which Monty Hall opened Door 2.

Notice that the game tree shows which actions Monty Hall and the contestant may take, but not which actions they \emph{should} take. For that, each player needs a strategy.

\begin{defn}
    A \emph{pure strategy} is a rule that tells a player how to move whenever it is his or her turn. 
    In other words, a pure strategy for player $p$ is a choice function
     	\[
    		\sigma \in \!\prod_{h \in H_{p}}\!A(h)
    	\]
    that respects the player's indistinguishability relation $\sim_{p}$, that is,
        \[
           	\sigma(h) = \sigma(h') \quad\text{whenever}\quad h \sim_{p} h'\!.
        \]
    The set of pure strategies for player $p$ is denoted $S_{p}$. 
\end{defn}

    A player $p$ is said to \emph{follow} a strategy \(\sigma \in S_{p}\) in a history \(h' \in H\) if, 
    for every history \(h \in H_{p}\) that is a proper initial segment of $h'$\!, 
    the history \(h\concat\sigma(h)\) is also an initial segment of $h'$\!. 
    In a win-lose game, a pure strategy is \emph{winning} if its player wins every terminal history in which he or she follows it.
    
The outcome of an extensive game is determined once every player has chosen a strategy. For example, suppose Monty Hall decides to place the prize behind Door 1 and to open the lowest-numbered door possible (i.e., not containing the prize or selected by the contestant). For her part, if the contestant selects Door 1, then switches to whichever door is not opened by Monty Hall, she will lose by switching to Door 3 after Monty Hall opens Door 2.

\begin{defn}
    A \emph{pure-strategy profile} is a vector 
    \(\boldsymbol{\sigma} = \langle\, \sigma_{p} \in S_{p} : p \in N\,\rangle\) 
    of pure strategies, one for each player.
\end{defn}

\begin{defn}
    Let 
    \(S = \prod_{p\in N} S_{p}\) 
    be the set of all pure-strategy profiles, 
    and let 
    \(h_{\boldsymbol{\sigma}} \in Z\) 
    denote the terminal history induced by the pure-strategy profile 
    \(\boldsymbol{\sigma} \in S\). 
    By abuse of notation, we will let
    \(u(\boldsymbol{\sigma}) = u(h_{\boldsymbol{\sigma}})\)
    denote the payoff vector received by the players when they follow the pure-strategy profile $\boldsymbol{\sigma}$\!.
\end{defn}

Although the contestant in the Monty Hall problem does not have a winning pure strategy, she does have an effective counterstrategy to each of Monty Hall's strategies. For example, if Monty Hall always places the prize behind Door 1, then the contestant can win by choosing Door 1 and sticking with it when offered the opportunity to switch. Similarly, if the contestant always chooses Door 1 and sticks with it, Monty Hall can counter by placing the prize behind Door 2 or Door 3. Thus, the losing player can always improve his or her payoff by changing strategies.

\begin{defn}
    Let \(\boldsymbol\sigma = \langle\,\sigma_{p} : p \in N\,\rangle\) be a pure-strategy profile. 
    A \emph{unilateral deviation} from $\boldsymbol\sigma$ by player $p$ 
    is a pure-strategy profile \(\boldsymbol\sigma' = \langle\,\sigma_{p}' : p \in N\,\rangle\)
    in which every player other than $p$ follows the same pure strategy as in $\boldsymbol\sigma$, 
    i.e., for all \(q \in N\setminus\{p\}\), 
    we have \(\sigma_{q} = \sigma_{q}'\).
    We will follow the standard convention of writing \(\boldsymbol\sigma' = \langle\sigma_{p}', \sigma_{-p}\rangle\).
    Such a deviation is \emph{profitable} for player $p$ if 
    \(u_{p}(\boldsymbol\sigma) < u_{p}(\boldsymbol\sigma')\).
\end{defn}

\begin{defn}
    A \emph{pure-strategy equilibrium} \(\boldsymbol\sigma^{*}\!\in S\) is a pure-strategy profile from which no player has a profitable deviation.
\end{defn}

As we observed above, the Monty Hall problem does not have a pure-strategy equilibrium since there are profitable deviations from every strategy profile. However, if we allow the players to vary their pure strategies from one play to the next, we can study the long-run effect of selecting each pure strategy with a given probability.

\begin{defn}
    A \emph{mixed strategy} for player $p$ is a probability distribution over $S_{p}$. 
    The set of mixed strategies for player $p$ is denoted $\Delta(S_{p})$. 
\end{defn}

When the players follow mixed strategies instead of pure strategies, 
the outcome of the game is no longer determined, but each outcome will occur with a certain probability.

\begin{defn}
    A \emph{mixed-strategy profile} is a vector
    \(\boldsymbol{\mu} = \bigl\langle\,\mu_{p} \in \Delta(S_{p}) : p \in N\,\bigr\rangle\)
    of mixed strategies, one for each player.
    If we assume that the players select their pure strategies independently, 
    then each mixed-strategy profile induces a probability distribution $\mu$ on  
    \(S = \prod_{p\in N} S_{p}\)
     defined by
	\[
		\mu(\boldsymbol{\sigma}) = 	\prod_{p\in N}\mu_{p}(\sigma_p),
	\]    
    where 
    \(\boldsymbol{\sigma} = \langle\,\sigma_{p} \in S_{p} : p \in N\,\rangle\).
    When $S$ is finite, we can define the \emph{expected utility} of $\boldsymbol\mu$ by
        \[
            	U(\boldsymbol{\mu}) = 
            	\sum_{\boldsymbol{\sigma}\in S} 
            	\mu(\boldsymbol{\sigma})
            	u(\boldsymbol{\sigma}),
        \]
    where the \emph{expected utility for player} $p$ is 
    \(U_{p}(\boldsymbol{\mu}) = U(\boldsymbol{\mu})(p)\). 
\end{defn}

\begin{example}
    Consider a two-player, win-lose game in which the first player has three pure strategies, 
    \(S_{I} = \{\sigma_{1}, \sigma_{2}, \sigma_{3}\}\), 
    and the second player has only two pure strategies,
    \(S_{II} = \{\tau_{1}, \tau_{2}\}\). 
    Suppose the first player follows a mixed strategy $\mu$ such that 
    	\(
    		\mu(\sigma_{1}) = \mu(\sigma_{2}) = \mu(\sigma_{3}) = 1/3,
    	\) 
    while the second player follows a mixed strategy $\nu$ such that 
    \(\nu(\tau_{1}) = 2/5\), 
    and 
    \(\nu(\tau_{2}) = 3/5\).
    Figure \ref{Figure-Mixed Strategy} depicts the situation where
    \begin{align*}
    	u_{I}(\sigma_{1},\tau_{1}) = 
    	u_{I}(\sigma_{2},\tau_{2}) = 
    	u_{I}(\sigma_{3},\tau_{1}) &= 
    	1, \\
    	u_{I}(\sigma_{1},\tau_{2}) = 
    	u_{I}(\sigma_{2},\tau_{1}) = 
    	u_{I}(\sigma_{3},\tau_{2}) &= 
    	0. 
    \end{align*}
    The area of the shaded region is equal to 
    \(U_{I}(\mu,\nu) = 7/15\), 
    while the area of the unshaded region is 
    \(U_{II}(\mu,\nu) = 8/15\). 
    
    One can easily see that neither player's mixed strategy is optimal. 
    If player II follows $\nu$, 
    then player I can improve her chances of winning by following $\sigma_{2}$ more often and 
    $\sigma_{1}$ and $\sigma_{3}$ less often. 
    Conversely, if player I follows $\mu$, 
    then player II can improve his expected utility  by following 
    $\tau_{2}$ more often and 
    $\tau_{1}$ less often.
\end{example}

\begin{figure}[htbp]
	\centering
	\begin{tikzpicture} 
		\fill (0,0) -- (0,1.5) -- (2.0,1.5) -- (2.0,0) [fill = gray!50];
		\fill (2.0,1.5) -- (2.0,3.0) -- (4.5,3.0) -- (4.5,1.5) [fill = gray!50];
		\fill (0,3.0) -- (0,4.5) -- (2.0,4.5) -- (2.0,3.0) [fill = gray!50];
		\draw (0,0) rectangle (4.5,4.5);
		\draw (2.0,0) -- (2.0,4.5); 
		\draw (0,1.5) -- (4.5,1.5); 
		\draw (0,3.0) -- (4.5,3.0); 
		%
		\draw (-0.4,3.75) node {$\sigma_{1}$};
		\draw (-0.4,2.25) node {$\sigma_{2}$};
		\draw (-0.4,0.75) node {$\sigma_{3}$};
		\draw (1.0,4.9) node {$\tau_{1}$};
		\draw (3.25,4.9) node {$\tau_{2}$};
	\end{tikzpicture} 
	\caption{A pair of mixed strategies}
  \label{Figure-Mixed Strategy}
\end{figure} 		


The previous example shows that players can profitably deviate from mixed-strategies as well as pure-strategies.

\begin{defn}
    Let \(\boldsymbol\mu = \bigl\langle\,\mu_{p} \in \Delta(S_{p}) : p \in N\,\bigr\rangle\) be a mixed-strategy profile. 
    A \emph{unilateral deviation} from $\boldsymbol\mu$ by player $p$ 
    is a mixed-strategy profile \(\boldsymbol\mu' = \bigl\langle\,\mu_{p}' \in \Delta(S_{p}) : p \in N\,\bigr\rangle\)
    in which every player other than $p$ follows the same mixed strategy as in $\boldsymbol\mu$, 
    i.e., for all \(q \in N\setminus\{p\}\), 
    we have \(\mu_{q} = \mu_{q}'\).
    We will follow the standard convention of writing \(\boldsymbol\mu' = \langle\mu_{p}', \mu_{-p}\rangle\).
    Such a deviation is \emph{profitable} for player $p$ if 
    \(U_{p}(\boldsymbol\mu) < U_{p}(\boldsymbol\mu')\).
\end{defn}

\begin{defn}
    A \emph{mixed-strategy equilibrium} is a mixed-strategy profile \(\boldsymbol\mu \in \prod_{p \in N} \Delta(S_{p})\) 
    from which no player has a profitable deviation. 
\end{defn}

John Nash proved that a strategic game has a mixed-strategy equilibrium if there are finitely many players that each have a finite number of pure strategies \cite{Nash:1950b,Nash:1951}. Nash's theorem is a generalization of von Neumann's minimax theorem, which states that every two-player, constant-sum game in which each player has a finite number of pure strategies has a mixed-strategy profile that simultaneously maximizes the minimum utility expected by each player. Player I's expected utility from such a mixed-strategy equilibrium is called the minimax value of the game \cite{Neumann:1928}. 

We now verify that the minimax value of the Monty Hall problem is 2/3. Observe that Monty Hall has a total of  \(3\cdot 2^{3} = 24\) pure strategies. However, pure strategies that differ only at decision points that are never reached when  those strategies are followed are \emph{outcome equivalent} \cite[page 94]{Osborne:1994}. 
By identifying outcome-equivalent strategies, we can reduce the number of Monty Hall's strategies by a factor of four. Let $\tau_{a}$ denote the \emph{reduced strategy} \cite[page 94]{Osborne:1994} according to which Monty Hall places the prize behind Door $a$, then opens the door with the lowest number possible, and let $\tau^{a}$ be the reduced strategy according to which he places the prize behind Door $a$, then opens the door with the highest number possible. 
Let $\nu^{*}$ be the mixed strategy according to which 
\(\nu^{*}(\tau_{a}) = 1/6 = \nu^{*}(\tau^{a})\). 


The contestant has \(3\cdot2^{6} = 192\) pure strategies, but only twelve reduced strategies.
Let $\sigma_{b}$ denote the reduced strategy according to which the contestant initially chooses Door $b$, then sticks with her initial choice when she is offered the opportunity to switch doors. Let $\sigma_{b}'$ denote the reduced strategy according to which the contestant initially chooses Door $b$, then switches doors. 
Observe that the contestant has six additional reduced strategies (that will remain nameless) according to which her decision to stick or switch doors depends on the door opened by Monty Hall.
Finally, let $\mu^{*}$ denote the mixed strategy according to which 
\(\mu^{*}\bigl(\sigma_{b}'\bigr) = 1/3\).

If the contestant sticks with her original door, she will win whenever it contains the prize, while, if she switches doors, she will win whenever her original door does not contain the prize. This information is summarized in Figure \ref{Figure-MHPayoffMatrix}, from which it is straightforward to calculate her expected utility \(U_{I}(\mu^{*}\!, \nu^{*}) = 2/3\). To show that neither player has a profitable deviation in mixed strategies, it suffices to fix $\mu^{*}$ and compute how well it fares against each of Monty Hall's reduced strategies, then fix $\nu^{*}$ and compute how it fares against each of the contestant's reduced strategies \cite[page 151]{Maschler:2013}. Again, inspecting Figure \ref{Figure-MHPayoffMatrix} reveals that 
    \[
	\min_{\tau \in S_{\mathrm{II}}} U_{I}(\mu^{*}\!, \tau) 
	= \frac{2}{3} 
	= \max_{\sigma \in S_{\mathrm{I}}} U_{I}(\sigma, \nu^{*}).
    \]
(Note that for each of the reduced strategies $\sigma$ not listed in Figure \ref{Figure-MHPayoffMatrix} we have \(U_{I}(\sigma, \nu^{*}) = 1/2\).)
Thus $\langle\mu^{*}\!, \nu^{*}\rangle$ is a mixed-strategy equilibrium, and the contestant's minimax value is 2/3.%
\footnote{
A similar analysis of the Monty Hall problem as an extensive game appears in Sandu \cite[pages 239--244]{Sandu:2015}.}


\begin{figure}[htbp]
\begin{center}
\[
\begin{array}{c|cccccc}
      			& \tau_{1}  & \tau_{2}  & \tau_{3}  & \tau^{1}  & \tau^{2}  & \tau^{3}   \\\hline
    \sigma_{1}  	& 1  & 0  & 0  & 1  & 0  & 0   \\
    \sigma_{2}  	& 0  & 1  & 0  & 0  & 1  & 0   \\
    \sigma_{3}  	& 0  & 0  & 1  & 0  & 0  & 1   \\
    \sigma_{1}'  	& 0  & 1  & 1  & 0  & 1  & 1   \\
    \sigma_{2}'  	& 1  & 0  & 1  & 1  & 0  & 1   \\
    \sigma_{3}'  	& 1  & 1  & 0  & 1  & 1  & 0   
\end{array}
\]

\end{center}
	\caption{Half of the reduced strategic form of the Monty Hall problem. The six rows that are not shown each contain three 0's and three 1's.}
  \label{Figure-MHPayoffMatrix}
\end{figure}


\section{First-order logic with imperfect information} 
\label{Section:IFLogic}

In this section, 
we introduce the syntax and game-theoretic semantics 
of first-order logic with imperfect information. 
There are several variants found in the literature,\!%
\footnote{
    In particular, 
    \vaananen's dependence logic \cite{Vaananen:2007} 
    has attracted a significant following.
    In dependence logic,
    the dependence relation between quantified variables is specified by 
    \emph{dependence atoms} such as 
    \[
    	=\!(x_{1},\ldots,x_{n},y),
    \]
    which indicates that the value of $y$ is determined by the values of 
    \(x_{1},\ldots,x_{n}\). 
}
but we adopt the original slashed notation of independence-friendly (IF) logic 
used by Hintikka and Sandu \cite{Hintikka:1989,Hintikka:1996}.
Our presentation will necessarily be brief. For a fuller treatment, 
we refer the reader to \cite{Mann:2011}.

First-order logic with imperfect information 
is a conservative extension of first-order logic that includes formulas whose semantic games are extensive games with imperfect information. 
Ordinary first-order formulas are built up from atomic formulas using 
negation ($\neg$),
disjunction ($\lor$),
conjunction ($\land$), 
existential quantification ($\exists$),
and 
universal quantification ($\forall$).
Atomic formulas and negated atomic formulas are called \emph{literals}. 
Independence-friendly formulas are similar to first-order formulas 
except that each connective and quantifier is parameterized by a finite set of variables that specifies the information available to the relevant player. 

\begin{defn}
Given any first-order vocabulary, an 
\emph{independence-friendly formula}
is an element of the smallest set
$\mathrm{IF}$ 
satisfying the following conditions:
 \begin{itemize}
  %
	\item 
	Every first-order literal belongs to $\mathrm{IF}$.
	%
	%
	\item 
	If \(\varphi,\psi \in \mathrm{IF}\), 
	and $W$ is a finite set of variables,
	then
	$(\varphi \hor{W} \psi)$ 
	and
	$(\varphi \hand{W} \psi)$
	belong to
	$\mathrm{IF}$.
	%
	\item 
	If \(\varphi \in \mathrm{IF}\), $x$ is a variable, 
	and $W$ is a finite set of variables, 
	then 
	\((\exists x/W)\varphi\) 
	and 
	\((\forall x/W)\varphi\)
	belong to 
	$\mathrm{IF}$.
  \end{itemize}
The finite set of variables $W$ 
is called an \emph{independence set} or \emph{slash set}.
\end{defn}
Informally,
a slash set indicates the variables in which 
a player's choice must be uniform. 
For example,
in the semantic game for 
\(\varphi \hor{\{x\}} \psi\),
\eloise's choice of disjunct may not depend on $x$.
In the semantic game for 
\(\bigl(\forall z/\{x,y\}\bigr)\varphi\), 
\abelard's choice of $z$ must be uniform in $x$ and $y$.
When a slash set is empty we simply omit it.
Free and bound variables are defined as usual, with the proviso that variables in slash sets are free.\!%
\footnote{
    For example,  
    the variable $x$ is free in the formula
    \(\bigl(\exists y/\{x\}\bigr)\varphi\),
    while both $x$ and $y$ are free in 
    \(\bigl(\exists y/\{x,y\}\bigr)\varphi\).
}
An \emph{independence-friendly sentence}
is an IF formula with no free variables.
The set of subformulas of an IF formula $\varphi$ is denoted 
$\subf(\varphi)$;
the set of literal subformulas of $\varphi$ is denoted 
$\lit(\varphi)$.

To avoid unnecessary complications, 
we will assume that all IF formulas are in negation normal form, i.e., the negation symbol only appears in front of atomic formulas.
However, for any IF formula $\varphi$,
we will use the notation $\neg\varphi$  
as a recursively defined abbreviation:
\begin{align*}
	\neg\neg\varphi & \quad\text{is}\quad \varphi, \\
	\neg(\varphi \hor{W} \psi) & \quad\text{is}\quad 
	\neg\varphi \hand{W} \neg\psi, \\
	\neg(\varphi \hand{W} \psi) & \quad\text{is}\quad 
	\neg\varphi \hor{W} \neg\psi, \\
	\neg(\exists x/W)\varphi & \quad\text{is}\quad 
	(\forall x/W)\neg\varphi, \\
	\neg(\forall x/W)\varphi & \quad\text{is}\quad 
	(\exists x/W)\neg\varphi. 
\end{align*} 
We will also use 
\(\varphi \implies \psi\)
as an abbreviation for 
\(\neg\varphi \lor \psi\).

Now that we have defined the syntax of IF logic, we next present its game-theoretic semantics.

\begin{defn} \label{Def:IFSemanticGame}
Let $\varphi$ be an IF sentence, and 
let $\str{M}$ be a suitable structure.\!%
\footnote{
    A structure is \emph{suitable} for an IF formula $\varphi$ 
    if it interprets every 
    function, 
    relation, and 
    constant symbol appearing in $\varphi$.
}
The \emph{semantic game} 
$G(\str{M},\varphi)$ is defined as follows:
\begin{itemize}
	\item 
	There are two players, 
	\eloise\ ($\exists$) 
	and
	\abelard\ ($\forall$).

\item The set of histories is $H = \bigcup \bigl\{\,H_{\psi } : \psi \in \subf(\varphi)\,\bigr\}$, where $H_{\psi }$ is defined recursively:
  \begin{itemize} 
  %
	\item $H_{\varphi} = \bigl\{(\varphi)\bigr\}$.
	%
	%
	\item If $\psi $ is $\chi_1 \hor{W} \chi_2$, then 
		\(
			H_{\chi_i }=\{\,h\concat\chi_i : h\in H_{\chi_1 \hor{W} \chi_2}\,\}
		\).
	%
	\item If $\psi $ is $\chi_1 \hand{W} \chi_2$, then 
		\(
			H_{\chi_i }=\{\,h\concat\chi_i : h\in H_{\chi_1 \hand{W} \chi_2}\,\}
		\).
	\item If $\psi $ is $(\exists x/W)\chi$, then 
	  \(
			H_{\chi } = 
			\bigl\{\,h\concat(x,a) : h\in H_{(\exists x/W)\chi},\, a\in M\,\bigr\}
	  \).
	\item If $\psi $ is $(\forall x/W)\chi$, then 
	  \(
			H_{\chi } = 
			\bigl\{\,h\concat(x,a) : h\in H_{(\forall x/W)\chi},\, a\in M\,\bigr\}
	  \).
\end{itemize}
Every history $h$ induces an assignment $s_{h}$ defined by
	\[
		s_{h}=
		\begin{cases}
    \emptyset			& \text{if $h = (\varphi)$}, \\
    s_{h'}  				& \text{if $h = h'\concat\psi$}, \\
    s_{h'}(x/a)  	& \text{if $h = h'\concat(x,a)$},
		\end{cases}
	\]
where $s_{h'}(x/a)$ is the assignment that is identical to $s_{h'}$ except that 
it assigns the value $a$ to the variable $x$. 
For example, 
let  $R(x,y)$ be an atomic formula, 
and
suppose $\varphi$ is
    \[
	\forall x\bigl(\exists y/\{x\}\bigr) \bigl[R(x,y) \land \neg R(x,y)\bigr].
    \] 
For any \(a,b\in M\),
the sequence 
\(h = \bigl(\varphi, (x,a), (y,b), \neg R(x,y)\bigr)\)
is a history for 
$G(\str{M},\varphi)$
that induces the assignment defined by 
    \[
	s_{h}(x) = a \quad\text{and}\quad
	s_{h}(y) = b.
    \]
	
	\item 
	Once play reaches a literal the game ends, 
	i.e., the set of terminal histories is:
		\[
			Z =\!\!\!\bigcup_{\chi \in \lit(\varphi)}\!\!\!\! H_{\chi}.
		\]
Observe that the above history $h$ is terminal because $\neg R(x,y)$ is a literal.

	\item 
	Disjunctions and existential quantifiers are decision points for 
	\eloise, 
	while conjunctions and universal quantifiers are decision points for 
	\abelard:
	\[
		P(h)=
		\begin{cases}
  	  \exists  
	  	&	\text{if 	$h\in H_{\chi_{1} \lor_{\!/W} \chi_{2}}$ or 
									$h\in H_{(\exists x/W)\chi}$}, \\
  	  \forall  
	  	& \text{if 	$h\in H_{\chi_{1} \land_{/W} \chi_{2}}$ or 
									$h\in H_{(\forall x/W)\chi }$}.
		\end{cases}
	\]
Let 
\(H_{\exists} = P^{-1}(\exists)\) 
and
\(H_{\forall} = P^{-1}(\forall)\)
be the sets of histories in which \eloise\ and \abelard\ are respectively active.
	\item 
	The indistinguishability relations 
	$\sim_{\exists}$ 
	and 
	$\sim_{\forall}$	
	are defined as follows. 
	\begin{itemize}
		\item
		For  
		\(h,h' \in H_{\chi_{1} \hor{W} \chi_{2}}\) or 
		\(h,h' \in H_{(\exists x/W)\chi}\), we have 
		\(h \sim_{\exists} h'\)
		if and only if
		\[
			\bigl\{\,x \in \dom(s_{h}) : s_{h}(x) \not= s_{h'}(x)\,\bigr\} \subseteq W.
		\]

		\item
		For  
		\(h,h' \in H_{\chi_{1} \hand{W} \chi_{2}}\) or 
		\(h,h' \in H_{(\forall x/W)\chi}\), we have 
		\(h \sim_{\forall} h'\) 
		if and only if
		\[
			\bigl\{\,x \in \dom(s_{h}) : s_{h}(x) \not= s_{h'}(x)\,\bigr\} \subseteq W.
		\]
	\end{itemize}
	Returning to our example $\phi$ above, 
	if \(a,a' \in M\),  
	then the histories 
	\(
		\bigl(\varphi, (x,a)\bigr) \sim_{\exists} \bigl(\varphi, (x,a')\bigr)
	\) 
	are indistinguishable to \eloise.
	
	\item 
	\eloise\ wins a terminal history \(h \in H_{\chi}\) 
	if \(\str{M}, s_{h} \models \chi\); 
	\abelard\ wins if \(\str{M}, s_{h} \not\models \chi\).
	
	For example,  
	\eloise\ wins the above history $h$ if and only if
	\((a,b) \notin R^{\str{M}}\)\!. 
\end{itemize}
\end{defn}

As with first-order logic, 
the truth or falsity of an IF sentence is determined not by winning a single play of the game, 
but by having a winning strategy.

\begin{defn} \label{Def:GTS}
Let $\varphi$ be an IF formula, and let $\str{M}$ be a suitable structure.
Then $\phi$ is \emph{true} in $\str{M}$, 
denoted \(\str{M} \modelt\!\phi\), 
if \eloise\ has a winning strategy for $G(\str{M},\phi)$,
and it is \emph{false} in $\str{M}$, 
denoted \(\str{M} \modelf\!\phi\), 
if \abelard\ has a winning strategy.
\end{defn}

The semantic game for a first-order sentence is a two-player, win-lose game with perfect information and finite horizon. Thus, the principle of bivalence for first-order logic is a consequence of the Gale--Stewart theorem \cite{Gale:1953}. Since the semantic game for an IF sentence is a game with imperfect information, the Gale--Stewart theorem does not apply. Hence, it is possible to have IF sentences that are neither true nor false in a given structure.

\begin{example}
    Consider the semantic game for the first-order sentence 
    \(
    	\forall x\exists y (x = y).
    \)
    When played on the two-element structure
    \(\str{2} = \{0,1\}\),
    \abelard\ has two possible strategies, 
    \(x \colonequals 0\) and 
    \(x \colonequals 1\),
    while \eloise\ has four possible strategies,
    \(y \colonequals 0\), 
    \(y \colonequals 1\),
    \(y \colonequals x\), 
    and 
    \(y \colonequals 1 - x\).
    Here, \(x \colonequals 0\) denotes the pure strategy that assigns the value 0 to $x$, 
    while $y \colonequals x$ denotes the pure strategy that assigns $y$ the same value as $x$.
    
    In contrast, \eloise\ cannot distinguish the histories 
    \(
    	\bigl(\phi, (x,0)\bigr) \sim_{\exists}
    	\bigl(\phi, (x,1)\bigr)
    \)
in the semantic game for the IF sentence 
    \(
    	\forall x\bigl(\exists y/\{x\}\bigr)\,x = y.
    \)
    Thus she only has two strategies, 
    \(y \colonequals 0\) and 
    \(y \colonequals 1\),
    neither of which is winning.
    It is easy to see that neither of Abelard's strategies, 
    \(x \colonequals 0\) and 
    \(x \colonequals 1\),
    are winning either
    (see Figure \ref{Figure-NoWinningStrategy}).

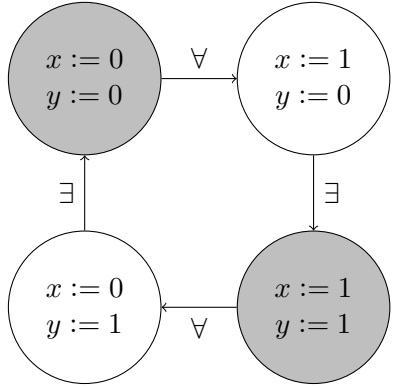
\begin{figure}[htbp]
	\centering
		\begin{tikzpicture}
			\node (00)	[circle,fill=gray!50,draw] {
				\(
					\begin{array}{c}
						x \colonequals 0 \\
						y \colonequals 0
					\end{array}
				\)};
			\node (10)			[right=of 00]	 [circle,draw] {				
				\(
					\begin{array}{c}
						x \colonequals 1 \\
						y \colonequals 0
					\end{array}
				\)}
				edge [<-] node[auto, swap] {$\forall$} (00);
			\node (11)			[below=of 10]	 [circle,fill=gray!50,draw] {				
				\(
					\begin{array}{c}
						x \colonequals 1 \\
						y \colonequals 1
					\end{array}
				\)}
				edge [<-] node[auto, swap] {$\exists$} (10);
			\node (01)			[left=of 11]	[circle,draw] {				
				\(
					\begin{array}{c}
						x \colonequals 0 \\
						y \colonequals 1
					\end{array}
				\)}
				edge [<-] node[auto, swap] {$\forall$} (11)
				edge [->] node[auto] {$\exists$} (00);

		\end{tikzpicture}
			\caption{The strategic form of the semantic game for 
				\(\forall x\bigl(\exists y/\{x\}\bigr)\,x = y\) when played on the structure \(\str{2} = \{0,1\}\). The arrows represent profitable deviations for the indicated player.}
  \label{Figure-NoWinningStrategy}
\end{figure}
\end{example}
 
Notice that the strategic form of the semantic game in the previous example is equivalent to the game Matching Pennies, in which two players simultaneously turn their respective coins to Heads or Tails. The first player wins if the coins match;
the second player wins if they differ. Although Matching Pennies does not have a pure-strategy equilibrium, it does have a mixed-strategy equilibrium where both players turn their coins to Heads or Tails with equal probability. 
The semantic game for \(\forall x\bigr(\exists y/\{x\}\bigr)\,x = y\) 
played on the structure \(\str{2} = \{0,1\}\) 
has a similar mixed-strategy equilibrium $\langle\mu^{*}\!, \nu^{*}\rangle$, 
where 
\begin{align*}
    \nu^{*}(x \colonequals 0) &= 1/2 = \nu^{*}(x \colonequals 1), \\
    \mu^{*}(y \colonequals 0) &= 1/2 = \mu^{*}(y \colonequals 1),
\end{align*}
and the minimax value for \eloise\ is \(U_{\exists}(\mu^{*}\!, \nu^{*}) = 1/2\). 
In a structure with $n$ elements, 
the minimax value for \eloise\ is $1/n$.%
\footnote{
    This example, due to Ajtai, first appeared in  
    \cite{Blass:1986}.
}

\begin{defn}
    The \emph{truth value} of an IF sentence $\varphi$ 
    in a suitable finite structure $\str{M}$ 
    is the minimax value for \eloise\ 
    in the semantic game $G(\str{M},\varphi)$.
\end{defn}

The above definition is due to Sevenster and Sandu \cite{Sevenster:2010}, who dub their extension of the basic game-theoretic semantics for IF logic \emph{equilibrium semantics}. Galliani \cite{Galliani:2008} independently developed a similar semantics based on behavioral strategies.%
\footnote{
    A \emph{behavioral strategy} for player $p$ 
    is a function mapping each of his or her information sets 
    to a probability distribution over the set of possible actions 
    at that information set.
}

We are now ready to show that, 
when interpreted in the structure
\(\str{M} = \{1,2,3\}\), 
the semantic game for the sentence
\(\varphi_{\mathrm{MH}}\),
\[
	\forall x 
	\bigl(\exists y/\{x\}\bigr)
	\forall z
	\Bigl[\,
		(z \not= x \land z \not= y) 
		\implies 
		\bigl(\exists y/\{x\}\bigr)\,
		x = y\,
	\Bigr], 
\]
is equivalent to the Monty Hall problem. 
Note that the above sentence is an abbreviation for
\[
	\forall x 
	\bigl(\exists y/\{x\}\bigr)
	\forall z
	\Bigl[\,
		(z = x \lor z = y)
		\lor 
		\bigl(\exists y/\{x\}\bigr)\,
		x = y\,
	\Bigr].
\]
At first glance, 
the semantic game 
$G(\str{M},\varphi_{\mathrm{MH}})$ 
appears more complicated than the Monty Hall problem 
since the contestant only makes two moves,
whereas \eloise\ makes three moves.
Moreover, 
\abelard\ and \eloise\ each have more possible actions in 
$G(\str{M},\varphi_{\mathrm{MH}})$ 
than Monty Hall and the contestant, respectively.  
Figure \ref{Figure-PhiMHGameTree}
shows the part of 
$G(\str{M},\varphi_{\mathrm{MH}})$ 
that occurs after \abelard\ sets the value of $x$ to 1, 
which corresponds to the upper third of 
Figure \ref{Figure-MontyHallGameTree}.

\input{Figure-PhiMHGameTree}

\abelard's possible strategies in the game
$G(\str{M},\varphi_{\mathrm{MH}})$ 
are straightforward.
First, he picks one of three possible values for $x$.
Then, after \eloise\ chooses a value for $y$,
he picks a value for $z$
(that may depend on $x$ and $y$).
Thus, a pure strategy for \abelard\
can be encoded by a pair $(a,f)$, 
where 
\(a \in \{1,2,3\}\)
and
\(f\colon \{1,2,3\}^{2} \to \{1,2,3\}\), 
while every such pair corresponds to a different pure strategy.
It follows that \abelard\ has  
\(3\cdot 3^{3^{2}} = 59,\!049\)
pure strategies.
However, 
suppose \abelard\ follows the pure strategy encoded by $(a,f)$.
Then, for all \(a' \not= a\), 
the value of $f(a'\!,b)$
is irrelevant because the history 
\(
	\bigl(\varphi_{\mathrm{MH}}, (x,a'), (y,b)\bigr)
\)
will never be played.
Thus we may assume that, 
for all \(a'\!,b \in \{1,2,3\}\),
we have \(f(a'\!,b) = f(a,b)\), 
reducing the number of \abelard's pure strategies to 
\(3\cdot 3^{3} = 81\).

Next observe that, if \abelard\ assigns $z$ the same value as $x$ or $y$, 
then \eloise\ can win by choosing the appropriate disjunct.
Consequently, 
\abelard\ should only follow strategies that lead 
to histories in which the value of $z$ is distinct from the values of $x$ and $y$,
of which there are only six.
Let $\tau_{a}$ and $\tau^{a}$ denote the strategies defined by
\(\tau_{a}(\varphi_{\mathrm{MH}}) = (x,a) = \tau^{a}(\varphi_{\mathrm{MH}})\) and 
\begin{align*}
	\tau_{a}\big(\varphi_{\mathrm{MH}}, (x,a), (y,b)\bigr) &= 
	\Bigr(z,\min\bigl(\{1, 2, 3\}\setminus \{a,b\}\bigr)\Bigr), 
	\\
	\tau^{a}\big(\varphi_{\mathrm{MH}}, (x,a), (y,b)\bigr) &= 
	\Bigr(z,\max\bigl(\{1, 2, 3\}\setminus \{a,b\}\bigr)\Bigr). 
\end{align*}

\eloise\ has even more strategies than \abelard.
Each of \eloise's pure strategies can be encoded by a quadruple $(b,g,h,j)$, where 
\(b \in \{1,2,3\}\) indicates the initial value she assigns to $y$,
the functions
\begin{align*}
	g&\colon \{1,2,3\}^{3} \to 
	\bigl\{(x = z\lor y = z), (\exists y/x)\,x = y\bigr\} \\
	h&\colon \{1,2,3\}^{3} \to 
	\bigl\{x = z, y = z\bigr\}
\end{align*}
indicate her choices of disjuncts,
and the function
\[
	j\colon \{1,2,3\}^{2} \to \{1,2,3\}
\]
indicates the final value she assigns to $y$.
Thus she has a total of 
\(3\cdot 2^{3^{3}}\!\cdot 2^{3^{2}}\!\cdot 3^{3^{2}}\!= 2^{36}\!\cdot 3^{10}\)
pure strategies.
Fortunately, 
we only need to consider a small fraction of them.
For starters, 
when faced with the disjunction
\[
		(x = z \lor y = z) 
		\lor 
		\bigl(\exists y/\{x\}\bigr)\,x = y\
\]
\eloise\ should only choose the left disjunct if 
the value of $z$ matches the value of $x$ or $y$.
Otherwise, she should choose the right disjunct. 
Hence 
there is only one function worth considering: 
\[
	g(a,b,c) = 
	\begin{cases}
  	(x = z \lor y = z)			& \text{if \(a = c\) or \(b = c\)}, \\
   	\bigl(\exists y/\{x\}\bigr)\,x = y		& \text{otherwise},
\end{cases}
\]
where 
$a$ is the value assigned to $x$, 
$b$ is the value assigned to $y$, and 
$c$ is the value assigned to $z$.
Thus, if \eloise\ ever finds herself faced with the disjunction $(x = z \lor y = z)$, 
the current assignment will satisfy one of the disjuncts.
Therefore, it is sufficient for her to always use the function 
\[
h(a,b,c) = 
	\begin{cases}
  	x = z										& \text{if \(a = c\)}, \\
  	y = z										& \text{otherwise}.
	\end{cases}
\]

When faced with the existential subformula \(\bigl(\exists y/\{x\}\bigr)\,x = y\),
\eloise\ should never set the final value of $y$ 
equal to the value of $z$ 
because she knows the values of $x$ and $z$ are distinct.
(Were they the same, 
she would have chosen the other disjunct.)
Similarly, 
she also knows that the values of $y$ and $z$ will differ. 
Thus 
she need only consider functions 
\[
	j\colon \Bigl(\{1,2,3\}^{2} \setminus \bigl\{(b,b) : b \in \{1,2,3\}\bigr\}\Bigr) \to \{1,2,3\}
\] 
such that \(j(b,c) \not= c\).
There are $2^{6}\!= 64$ such functions, 
two of which are of particular interest.
Let $j_{\mathrm{stick}}$ be the function such that for all \(b,c \in \{1,2,3\}\) we have
\(j_{\mathrm{stick}}(b,c) = b\),
and let $j_{\mathrm{switch}}$ be the function such that 
$j_{\mathrm{switch}}(b,c)$ 
is the unique element in \(\{1,2,3\}\setminus\{b,c\}\).

Let 
\(\sigma_{b}\) 
denote the pure strategy for \eloise\ encoded by the quadruple $(b,g,h,j_{\mathrm{stick}})$, 
and let 
\(\sigma_{b}'\)
be the strategy encoded by $(b,g,h,j_{\mathrm{switch}})$. 
Suppose that $\sigma$ and $\tau$ are rational strategies for \eloise\ and \abelard, respectively, and that the pair $(\sigma,\tau)$ induce the terminal history 
\[
	\Bigl(\varphi_{\mathrm{MH}}, 
	(x,a), 
	(y,b), 
	(z,c), 
	\bigl(\exists y/\{x\}\bigr)\,x = y, 
	(y,d)
	\Bigr),
\]
where \(a \not= c \not= b\) and \(c \not= d\). 
Observe that this is the same as the terminal history induced by 
$(\sigma_{b},\tau)$ or $(\sigma_{b}',\tau)$, 
depending on whether \(b = d\) or \(b \not= d\), respectively.

Figure \ref{Figure-PhiMHGameTree} highlights
the terminal histories that are induced by pairs of rational strategies.
By inspection, 
one can see that these histories form a subtree that is isomorphic 
to the game tree for the Monty Hall problem shown in 
Figure \ref{Figure-MontyHallGameTree}.
Moreover,
one can check that two histories are indistinguishable to the contestant
if and only the corresponding histories are indistinguishable to \eloise.
Thus, the Monty Hall problem has the same strategic form as the semantic game
$G(\str{M},\varphi_{\mathrm{MH}})$ 
if we assume that \abelard\ and \eloise\ follow rational strategies.
Every finite game has a mixed-strategy equilibrium involving only rational strategies 
\cite[Proposition 122.1]{Osborne:2004}. 
Thus, our previous analysis of the Monty Hall problem shows that 
the truth value of $\varphi_{\mathrm{MH}}$ is 2/3.%
\footnote{
A similar analysis of the Monty Hall problem as the semantic game of an IF sentence appears in Sandu \cite[pages 244--245]{Sandu:2015}.}

To conclude this section,
let us briefly analyze a variant of the Monty Hall problem in which the host is not required to offer the contestant the opportunity to switch doors. That is, Monty Hall is allowed to open any of the three doors, including the door initially chosen by the contestant or the door containing the prize.
\begin{enumerate}
    \item
    If the host opens the contestant's initial door or the door containing the prize, the contestant wins if she guessed correctly and loses if she did not. 
    \item
    If the host opens neither the contestant's initial door nor the door containing the prize, he then offers the contestant the opportunity to switch doors.
\end{enumerate}
This scenario can be modeled by the IF sentence 
$\varphi_{\mathrm{MH}}'$,
    \[
	\forall x 
	\bigl(\exists y/\{x\}\bigr)
	\forall z
	\biggl[
		(x = y = z) \lor
	\Bigl[\, 
		\neg(z = x \not = y) 
		\land 
		\neg(z = y \not= x) 
		\land 
		\bigl(\exists y/\{x\}\bigr)\,
		x = y\,
	\Bigr]
	\biggr]. 
    \]

\begin{figure}[htbp]
	\centering
		\resizebox{4in}{!}{
		\begin{tikzpicture}
			\node (root) 	{};
			\node (b)			[rounded corners,draw] 
										[right=of root]	{$\exists y/\{x\}$}
											edge[very thick] node[auto,swap] {1} (root);
			\node (b1)		[circle,draw] 
										[above right=of b, xshift=5mm, yshift=5mm] {$\forall z$}
											edge[very thick] node[auto,swap] {1} (b);
			\node (b11)	[circle,fill=gray]
										[above right=of b1, yshift=5mm, label=right:{$\exists$}] {}
											edge node[auto,swap] {1} (b1);
			\node (b12)	[circle,draw]
										[right=of b1, xshift=5mm] {$\lor$}
											edge[very thick] node[auto,swap] {2} (b1);
			\node (b12A)	[circle,fill=gray]
										[above right=of b12, yshift=0mm, label=right:{$\forall$}]{}
											edge node[auto,swap] {} (b12);
			\node (b12and)	[circle,draw]
										[right=of b12, xshift=5mm]{$\land$}
											edge[very thick] node[auto,swap] {} (b12);
			\node (b12andx)	[circle,fill=gray]
										[above right=of b12and, yshift=0mm, label=right:{$\exists$}]{}
											edge node[auto,swap] {} (b12and);
			\node (b12andy)	[circle,fill=gray]
										[right=of b12and, yshift=0mm, label=right:{$\exists$}]{}
											edge node[auto,swap] {} (b12and);
			\node (b12Ey)	[rounded corners,draw]
										[below right=of b12and]{$\exists y/\{x\}$}
											edge[very thick] node[auto,swap] {} (b12and);
			\node (b12Ey1)		[circle,fill] 
										[above right=of b12Ey, label=right:{$\exists$}]	{}
											edge[very thick] node[auto,swap] {1} (b12Ey);
			\node (b12Ey2)		[circle,fill] 
										[right=of b12Ey, label=right:{$\forall$}]	{}
											edge[very thick] node[auto,swap] {2} (b12Ey);
			\node (b12Ey3)		[circle,fill] 
										[below right=of b12Ey, label=right:{$\forall$}] {}
											edge[very thick] node[auto] {3} (b12Ey);
			\node (b13)	[circle,draw]
										[below right=of b1, yshift=-5,,] {$\lor$}
											edge[very thick] node[auto] {3} (b1);
			\node (b13A)	[circle,fill=gray]
										[right=of b13, yshift=0mm, label=right:{$\forall$}]{}
											edge node[auto,swap] {} (b13);
			\node (b13and)	[circle,draw]
										[below right=of b13, yshift=-4mm]{$\land$}
											edge[very thick] node[auto,swap] {} (b13);
			\node (b13andx)	[circle,fill=gray]
										[above right=of b13and, yshift=0mm, label=right:{$\exists$}]{}
											edge node[auto,swap] {} (b13and);
			\node (b13andy)	[circle,fill=gray]
										[right=of b13and, yshift=0mm, label=right:{$\exists$}]{}
											edge node[auto,swap] {} (b13and);
			\node (b13Ey)	[rounded corners,draw]
										[below right=of b13and, xshift=-4mm, yshift=-4mm]{$\exists y/\{x\}$}
											edge[very thick] node[auto,swap] {} (b13and);
			\node (b13Ey1)		[circle,fill] 
										[above right=of b13Ey, label=right:{$\exists$}]	{}
											edge[very thick] node[auto,swap] {1} (b13Ey);
			\node (b13Ey2)		[circle,fill] 
										[right=of b13Ey, label=right:{$\forall$}]	{}
											edge[very thick] node[auto,swap] {2} (b13Ey);
			\node (b13Ey3)		[circle,fill] 
										[below right=of b13Ey, label=right:{$\forall$}] {}
											edge[very thick] node[auto] {3} (b13Ey);
			\node (bb)		[circle,fill] 
										[right=of b, label=right:{$\forall$}]	{}
											edge[very thick] node[auto,swap] {2} (b);

			\node (bc)		[circle,fill] 
										[below right=of b, label=right:{$\forall$}] {}
											edge[very thick] node[auto] {3} (b);
		\end{tikzpicture}
		}
			\caption{A portion of the semantic game for \(\varphi_{\mathrm{MH}'}\)
			}
  \label{Figure-PhiMHprimeGameTree}
\end{figure}
 
 A portion of the game tree for $G\bigl(\str{M}, \varphi_{\mathrm{MH}}'\bigl)$ is shown in Figure \ref{Figure-PhiMHprimeGameTree}. If the initial values of all three variables are equal, \eloise\ can win by choosing the disjunct \(x = y = z\). In contrast, if \eloise\ initially assigns $y$ a different value than $x$, then \abelard\ can win by setting the value of $z$ equal to the value of $x$ or $y$. Thus, if $x$ and $y$ are initially assigned the same value, \abelard\ should pick a different value for $z$, forcing \eloise\ to choose the right disjunct, after which \abelard\ is forced to pick the rightmost conjunct, giving \eloise\ the opportunity to assign $y$ a new value that may depend on the values of $y$ and $z$, but not $x$. 

\newpage
To help define pure strategies for both players, let 
\(h_{a} = \bigl(\varphi_{\mathrm{MH}}', (x,a)\bigr)\), 
\(h_{ab} = h_{a}\concat(y,b)\),
and
\(h_{abc} = h_{ab}\concat(z,c)\).
\abelard's pure strategies $\tau_{a}$ and $\tau^{a}$ are defined as follows:
    \[
        \tau_{a}\bigl(\varphi_{\mathrm{MH}}'\bigr) = (x,a) = \tau^{a}\bigl(\varphi_{\mathrm{MH}}'\bigr),
    \]
    \begin{align*}
        \tau_{a}(h_{ab}) &= 
        \begin{cases}
            (z,a) & \text{if \(a \not= b\)}, \\
            \Bigr(z,\min\bigl(\{1, 2, 3\}\setminus \{a\}\bigr)\Bigr) & \text{otherwise},
        \end{cases} \\
        \tau^{a}(h_{ab}) &= 
        \begin{cases}
            (z,a) & \text{if \(a \not= b\)}, \\
            \Bigr(z,\max\bigl(\{1, 2, 3\}\setminus \{a\}\bigr)\Bigr) & \text{otherwise}, 
        \end{cases}
    \end{align*}
    \begin{align*}
        \tau_{a}\big(h_{abc},\ldots,\psi\big) = \tau^{a}\big(h_{abc},\ldots,\psi\big) &= 
        \begin{cases}
            \neg(z = x \not = y) & \text{if \(a = c\)}, \\
            \neg(z = y \not= x) & \text{if \(a \not= b = c\)}, \\
            \bigl(\exists y/\{x\}\bigr)x = y & \text{otherwise}, 
        \end{cases}
    \end{align*}
where $\psi$ is the conjunction 
    \[
        \neg(z = x \not = y) 
        \land 
        \neg(z = y \not= x) 
        \land 
        \bigl(\exists y/\{x\}\bigr)x = y.
    \]
Equivalently, we could have defined \(\tau_{a}(h_{ab}) = (z,b) = \tau^{a}(h_{ab})\) when \(a \not= b\).

Define \eloise's pure strategy $\sigma_{b}$ as follows: 
    \begin{align*}
        \sigma_{b}(h_{a}) &= (y,b), \\
        \sigma_{b}(h_{abc}) &=
            \begin{cases}
                x = y = z & \text{if \(a = b = c\),} \\
                \psi & \text{otherwise,}
            \end{cases} \\
        \sigma_{b}\Bigl(h_{abc},\ldots,\bigl(\exists y/\{x\}\bigr)\,x = y\,\Bigr) &= 
        (y, b),
    \end{align*}
and define $\sigma_{b}'$ similarly except that 
    \[
	\sigma_{b}'\Bigl(h_{abc},\ldots,\bigl(\exists y/\{x\}\bigr)\,x = y\,\Bigr) 
	= (y,b'),
    \]
where $b'$ is the unique element in $\{1, 2, 3\}\setminus\{b, c\}$ when $b$ and $c$ are distinct, otherwise \(b' = \min\bigl(\{1, 2, 3\}\setminus\{b\}\bigr)\). Equivalently, we could have defined \(b' = \max\bigl(\{1, 2, 3\}\setminus\{b\}\bigr)\) when \(b = c\).

Observe that, if \abelard\ follows $\tau_{a}$ or $\tau^{a}$ and \eloise\ follows $\sigma_{b}$ or $\sigma_{b}'$, then whenever play reaches the formula $\bigl(\exists y/\{x\}\bigr)x = y$, the values $x$ and $y$ will be the same while the value of $z$ will be different. Hence, she should not change the value of $y$.

Now define the mixed strategies \(\mu, \mu' \in \Delta(S_{\exists})\) and \(\nu \in \Delta(S_{\forall})\) by 
    \[
	\nu(\tau_{a}) = 1/6 = \nu(\tau^{a}),
    \]
    \[
	\mu(\sigma_{b}) = 1/3 = \mu'\bigl(\sigma_{b}'\bigr).
    \]
Then \(U_{\exists}(\mu,\nu) = 1/3\) because \eloise\ wins if and only if her initial guess is correct. It is easy to verify that $\langle{\mu,\nu}\rangle$ is a mixed-strategy equilibrium, since \abelard\ cannot improve his expected utility by favoring certain values of $x$ over  others. Nor does it matter which value he assigns to $z$ when the initial values of $x$ and $y$ are the same.
For her part, \eloise\ will win if and only if her initial guess is correct, assuming \abelard\ follows $\nu$. Thus the value of the semantic game $G\bigl(\str{M}, \varphi_{\mathrm{MH}}')$ is 1/3 (see Figure \ref{Figure-PhiMHprimeEquilibrium}).
 
\begin{figure}[htbp]
\begin{center}
\begin{tikzpicture} 
	\draw (2.25,5.0) node {\(x = 1\)};
	\fill (0,3.0) -- (0,4.5) -- (4.5,4.5) -- (4.5,3.0) [fill = gray!50];
	\draw (0,0) rectangle (4.5,4.5);
	\draw (2.25,3.0) -- (2.25,4.5); 
	\draw (0,1.5) -- (4.5,1.5); 
	\draw (0,3.0) -- (4.5,3.0); 
	\draw (1.125,3.75) node  
	  {\(
	  	\begin{array}{c}
		  y = 1 \\
		  z = 2
	  	\end{array}
	  \)};
	\draw (3.375,3.75) node  
	  {\(
	  	\begin{array}{c}
		  y = 1 \\
		 	z = 3
	  	\end{array}
	  \)};
	\draw (2.25,2.25) node  
	  {\(
	  	\begin{array}{c}
		  y = 2 \\
		  z = 1
	  	\end{array}
	  \)};
	\draw (2.25,0.75) node  
	  {\(
	  	\begin{array}{c}
		  y = 3 \\
		  z = 1
	  	\end{array}
	  \)};
	  
	\draw (7.75,5.0) node {\(x = 2\)};
	\fill (5.5,1.5) -- (5.5,3.0) -- (10.0, 3.0) -- (10.0,1.5) [fill = gray!50];
	\draw (5.5,0) rectangle (10.0,4.5); 
	\draw (5.5,1.5) -- (10.0,1.5); 
	\draw (5.5,3.0) -- (10.0,3.0); 
	\draw (7.75,1.5) -- (7.75,3.0); 
	\draw (7.75, 3.75) node  
	  {\(
	  	\begin{array}{c}
		  y = 1 \\
		  z = 2
	  	\end{array}
	  \)};
	\draw (6.675, 2.25) node  
	  {\(
	  	\begin{array}{c}
		  y = 2 \\
		  z = 1
	  	\end{array}
	  \)};
	\draw (8.875, 2.25) node  
	  {\(
	  	\begin{array}{c}
		  y = 2 \\
		  z = 3
	  	\end{array}
	  \)};
	\draw (7.75, 0.75) node  
	  {\(
	  	\begin{array}{c}
		  y = 3 \\
		  z = 1
	  	\end{array}
	  \)};
	  
	\draw (13.25, 5.0) node {\(x = 3\)};
	\fill (11,0) -- (15.5,0) -- (15.5, 1.5) -- (11,1.5) [fill = gray!50];
	\draw (11, 0) rectangle (15.5, 4.5) ; 
	\draw (11, 1.5) -- (15.5, 1.5); 
	\draw (11, 3.0) -- (15.5, 3.0); 
	\draw (13.25, 0.0) -- (13.25, 1.5);
	\draw (13.25, 3.75) node  
	  {\(
	  	\begin{array}{c}
		  y = 1 \\
		  z = 3
	  	\end{array}
	  \)};
	\draw (13.25, 2.25) node  
	  {\(
	  	\begin{array}{c}
		  y = 2 \\
		  z = 3
	  	\end{array}
	  \)};
	\draw (12.125, 0.75) node  
	  {\(
	  	\begin{array}{c}
		  y = 3 \\
		  z = 1
	  	\end{array}
	  \)};
	\draw (14.375, 0.75) node  
	  {\(
	  	\begin{array}{c}
		  y = 3 \\
		  z = 2
	  	\end{array}
	  \)};
\end{tikzpicture} 
\end{center}
\caption{A mixed-strategy equilibrium for $\varphi_{\mathrm{MH}}'$}
\label{Figure-PhiMHprimeEquilibrium}
\end{figure}


\section{Stochastic IF logic}
\label{Section:Stochastic IF logic}

Another variant of the Monty Hall problem involves treating the host as a disinterested party rather than as a malevolent opponent. 
If Monty Hall places the prize behind each door with equal probability and opens each door equally often (whether or not it conceals the prize or was chosen by the contestant), then, on those occasions when Monty happens to open a door that neither contains the prize nor was chosen by the contestant, the prize will be found behind each of the remaining doors with equal probability \cite[pages 935--936]{Friedman:1998}.

To model the scenario just described, we must go beyond ordinary IF logic by adding chance moves to our semantic games. We imagine that such moves are taken by a third player (Nature) who is indifferent to the eventual outcome of the game. Chance moves will be indicated by a new connective $\times$ and quantifier $\stochastic x$.%
\footnote{
    Sandu introduces what he calls \emph{probabilistic quantifiers}, denoted $\mu x$, where $\mu$ is a probability distribution over a finite universe \cite[page 246]{Sandu:2015}. We have adapted the backward-\textsf{S} notation from Alexey Radul's blog post: \\ 
    \texttt{http://alexey.radul.name/ideas/2014/stochasticity-is-a-quantifier/}
}

\begin{defn}
Given any first-order vocabulary, a \emph{stochastic IF formula} is a member of the smallest set $\mathrm{IF}(\stochastic)$ that satisfies the following conditions:
 \begin{itemize}
  %
	\item 
	Every IF formula belongs to IF($\stochastic$).
	%
	%
	\item 
	If \(\varphi,\psi \in \mathrm{IF}(\stochastic)\), then 
	\((\varphi \times \psi) \in \mathrm{IF}(\stochastic)\).	
	%
	\item 
	If \(\varphi \in \mathrm{IF}(\stochastic)\), and $x$ is a variable, then 
	\(\stochastic x\varphi \in \mathrm{IF}(\stochastic)\).
\end{itemize}
An \emph{stochastic IF sentence} 
is a stochastic IF formula with no free variables.
\end{defn}

\begin{defn}
Let $\varphi$ be a stochastic IF sentence, and 
let $\str{M}$ be a suitable structure. 
The \emph{semantic game}
$G(\str{M}, \varphi)$
is defined as before with the following amendments:
\begin{itemize}
	\item
	There are three players, Nature ($\stochastic$), \eloise\ ($\exists$), and \abelard\ ($\forall$).

	\item 
	If $\psi$ is \(\chi_{1} \times \chi_{2}\), 
	then \(H_{\chi_i} = \{\,h\concat\chi_{i} : h \in H_{\chi_{1} \times \chi_{2}}\,\}\).
	
	\item 
	If $\psi$ is $\stochastic x\chi$, then 
	\(H_{\chi} = \{\,h\concat (x,a) : h \in H_{\scriptsize{\stochastic x\chi}}, a \in M\,\}\).
	
	\item The player function is redefined to make the new connectives and quantifiers moves for Nature:
	\[
		P(h)=
		\begin{cases}
  	  \stochastic  
	  	& \text{if 	$h\in H_{\chi_{1} \times \chi_{2}}$ or 
									$h\in H_{\scriptsize{\stochastic} x\chi }$}, \\
  	  \exists  
	  	&	\text{if 	$h\in H_{\chi_{1} \lor_{\!/W} \chi_{2}}$ or 
									$h\in H_{(\exists x/W)\chi}$}, \\
  	  \forall  
	  	& \text{if 	$h\in H_{\chi_{1} \land_{/W} \chi_{2}}$ or 
									$h\in H_{(\forall x/W)\chi }$}.
		\end{cases}
	\]
	Let 
	\(H_{\scriptsize{\stochastic}} = P^{-1}(\stochastic)\) 
	be the set of histories in which Nature is the active player.
	
	\item 
	Nature's indistinguishability relation
	$\sim_{\scriptsize\stochastic}$ 
	is the identity relation. 
	That is, all of Nature's information sets are singletons.
	
	\item
	Nature receives no utility regardless of the outcome of the game. 
	For every terminal history \(h \in Z\), 
	we have \(u_{\scriptsize\stochastic}(h) = 0\).
\end{itemize}
\end{defn}

Since Nature has no reason to prefer one action over another, 
we cannot reason endogenously about which actions Nature will take. 
Instead, we will assume that Nature follows a behavioral strategy that is fixed in advance and known to all of the other players. Thus, we will treat Nature's strategy as an exogenous parameter. 

In an extensive game with chance moves, a pure-strategy profile for the players other than Nature does not uniquely determine a terminal history; it determines a probability distribution over the set of terminal histories. Thus, although we cannot predict the exact payoffs the players will receive based only on their own actions, we can compute the expected utility for each player given any mixed/behavioral-strategy profile. 

\begin{defn}
    The \emph{expectiminimax value} of a two-player, constant-sum extensive game with chance moves, 
    relative to a fixed behavioral strategy $\lambda$ for Nature, is
    \begin{align*}
    	\min_{\nu\in\Delta(S_{II})}\
    	\max_{\mu\in\Delta(S_{I})} U_{I}(\lambda,\mu,\nu)\,=\,
    	U_{I}\bigl(\lambda,\mu^{*}\!,\nu^{*}\bigr)\,=\!
    	\max_{\mu\in\Delta(S_{I})}\
    	\min_{\nu\in\Delta(S_{II})} U_{I}(\lambda,\mu,\nu).
    \end{align*}
\end{defn}

\begin{defn}
    Let $\varphi$ be a stochastic IF sentence, 
    let $\str{M}$ be a suitable finite structure, 
    and let $\lambda$ be a behavioral strategy for Nature in the semantic game $G(\str{M}, \varphi)$. 
    The \emph{truth value} of $\varphi$ in $\str{M}$ (relative to $\lambda$) 
    is the expectiminimax value of $G(\str{M}, \varphi)$ 
    when Nature follows $\lambda$.
    We will use the notation \(\str{M} \models^{v}\!\varphi(\lambda)\) 
    to express the fact that $v$ is the truth value of $\varphi$ in $\str{M}$ (relative to $\lambda$).
\end{defn}

\begin{example}
    Consider the following stochastic generalization of the Matching Pennies sentence:%
    \footnote{
    Sandu considers a different stochastic variant of the Matching Pennies sentence \cite[page 246]{Sandu:2015}.}
        \[
            \forall x\bigl(\exists y/\{x\}\bigr)\stochastic z [x = y = z]
        \]
    When played on the two-element structure $\str{2} = \{0, 1\}\), 
    the semantic game for the above sentence models the scenario 
    in which two players each turn a coin to Heads or Tails, 
    but now there is a third coin that Nature tosses in secret. 
    The first player (\eloise) wins if all three coins match; otherwise she loses. 
    
    If we assume that Nature's coin is biased, 
    so that Nature follows the behavioral strategy defined by \(\lambda(0) = 1/3\) and \(\lambda(1) = 2/3\), 
    then the semantic game can be represented by the game tree shown in 
    Figure \ref{Figure-StochasticMatchingPenniesGameTree}, 
    where \abelard\ follows the mixed strategy defined by \(\nu_{q}(0) = q\) and \(\nu_{q}(1) = 1 - q\), 
    while \eloise\ follows the mixed strategy defined by \(\mu_{p}(0) = p\) and \(\mu_{p}(1) = 1 - p\).
    Then \eloise's expected utility is \(U_{\exists}(\lambda, \mu_{p}, \nu_{q}) = \frac{1}{3}pq + \frac{2}{3}(1 - p)(1 - q)\).
    Using the second-derivative test, 
    we can show that there is an equilibrium when \(p = 2/3\) and \(q = 2/3\).
    Hence the truth value of the  sentence is \(U_{\exists}(\lambda, \mu_{2/3}, \nu_{2/3}) = 2/9\).
    This equilibrium is depicted in Figure \ref{Figure-StochasticMatchingPenniesEquilibrium}. 
    The reader should imagine that the square on the left labeled \(z = 0\) has a vertical ``thickness'' of 1/3, 
    while the square on the right labeled \(z = 1\) has a thickness of 2/3, 
    which accounts for the fact that Nature's coin is biased.
\end{example}


\begin{figure}[htbp]
    \begin{center}
        \begin{tikzpicture}
        			\node (root)	
                			[circle,draw] {$\forall x$};
			
    			\node (0)		
            			[rounded corners,draw] [above right=of root, yshift=25mm] {$\exists y/\{x\}$}
    				edge node[auto,swap] {$q$} (root);
    			\node (00)		
            			[circle,draw] [above right=of 0, yshift=5mm] {$\stochastic z$}
    				edge node[auto,swap] {$p$} (0);
    			\node (000)		
            			[circle,fill] [above right=of 00, label=right:{$\exists$}] {}
    				edge node[auto,swap] {$\frac{1}{3}$} (00);
    			\node (001)		
            			[circle,fill] [below right=of 00, label=right:{$\forall$}] {}
    				edge node[auto] {$\frac{2}{3}$} (00);
    			\node (01)		
            			[circle,draw] [below right=of 0, yshift=-5mm] {$\stochastic z$}
    				edge node[auto] {$1 - p$} (0);
    			\node (010)		
            			[circle,fill] [above right=of 01, label=right:{$\forall$}] {}
    				edge node[auto,swap] {$\frac{1}{3}$} (01);
    			\node (011)		
            			[circle,fill] [below right=of 01, label=right:{$\forall$}] {}
    				edge node[auto] {$\frac{2}{3}$} (01);
 
    			\node (1)		
            			[rounded corners,draw] [below right=of root, yshift=-25mm] {$\exists y/\{x\}$}
    				edge node[auto] {$1 - q$} (root)
				edge [dotted] node {} (0);
    			\node (10)		
            			[circle,draw] [above right=of 1, yshift=5mm] {$\stochastic z$}
    				edge node[auto,swap] {$p$} (1);
    			\node (100)		
            			[circle,fill] [above right=of 10, label=right:{$\forall$}] {}
    				edge node[auto,swap] {$\frac{1}{3}$} (10);
    			\node (101)		
            			[circle,fill] [below right=of 10, label=right:{$\forall$}] {}
    				edge node[auto] {$\frac{2}{3}$} (10);
    			\node (11)		
            			[circle,draw] [below right=of 1, yshift=-5mm] {$\stochastic z$}
    				edge node[auto] {$1 - p$} (1);
    			\node (110)		
            			[circle,fill] [above right=of 11, label=right:{$\forall$}] {}
    				edge node[auto,swap] {$\frac{1}{3}$} (11);
    			\node (111)		
            			[circle,fill] [below right=of 11, label=right:{$\exists$}] {}
    				edge node[auto] {$\frac{2}{3}$} (11);
        \end{tikzpicture}
    \end{center}
    \caption{The semantic game for the stochastic matching pennies sentence}
    \label{Figure-StochasticMatchingPenniesGameTree}
\end{figure}


\begin{figure}[htbp]
\begin{center}
\begin{tikzpicture} 
	\draw (2.25,5.0) node {\(z = 0\)};
	\fill (0,1.5) -- (3.0,1.5) -- (3.0, 4.5) -- (0.0, 4.5) [fill = gray!30];
	\draw (0,0) rectangle (4.5,4.5);
	\draw (0,1.5) -- (4.5,1.5); 
	\draw (3.0,0) -- (3.0,4.5); 
	\draw (1.5, 3.0) node  
	  {\(
	  	\begin{array}{c}
		  x = 0 \\
		  y = 0
	  	\end{array}
	  \)};
	\draw (1.5, 0.75) node  
	  {\(
	  	\begin{array}{c}
		  x = 0 \\
		  y = 1
	  	\end{array}
	  \)};
	\draw (3.75, 3.0) node  
	  {\(
	  	\begin{array}{c}
		  x = 1 \\
		  y = 0
	  	\end{array}
	  \)};
	\draw (3.75, 0.75) node  
	  {\(
	  	\begin{array}{c}
		  x = 1 \\
		  y = 1
	  	\end{array}
	  \)};
	  
	\draw (7.75,5.0) node {\(z = 1\)};
	\fill (8.5, 0.0) -- (10.0, 0.0) -- (10.0, 1.5) -- (8.5, 1.5) [fill = gray!30];
	\draw (5.5,0) rectangle (10.0,4.5); 
	\draw (8.5,0) -- (8.5, 4.5); 
	\draw (5.5, 1.5) -- (10.0, 1.5); 
	\draw (7.0, 3.0) node  
	  {\(
	  	\begin{array}{c}
		  x = 0 \\
		  y = 0
	  	\end{array}
	  \)};
	\draw (7.0, 0.75) node  
	  {\(
	  	\begin{array}{c}
		  x = 0 \\
		  y = 1
	  	\end{array}
	  \)};
	\draw (9.25, 3.0) node  
	  {\(
	  	\begin{array}{c}
		  x = 1 \\
		  y = 0
	  	\end{array}
	  \)};
	\draw (9.25, 0.75) node  
	  {\(
	  	\begin{array}{c}
		  x = 1 \\
		  y = 1
	  	\end{array}
	  \)};
	  
\end{tikzpicture} 
\end{center}
\caption{An equilibrium for the stochastic Matching Pennies sentence}
\label{Figure-StochasticMatchingPenniesEquilibrium}
\end{figure}

\newpage

We are now ready to model the variant of the Monty Hall problem in which the host is indifferent to the outcome as the semantic game of a stochastic IF sentence. We will assume that the contestant wins if Monty Hall reveals that the door she initially chose contains the prize and loses if he reveals that her door does not contain the prize or that some other door does contain the prize. Otherwise, she is offered the opportunity to switch doors. We will further assume that the prize is placed behind each door with probability 1/3, and that Monty Hall opens each door with probability 1/3.

 We will denote the following stochastic variant of $\varphi_{\mathrm{MH}}'$ by $\varphi_{\mathrm{MH}\scriptsize{\stochastic}}'$,
    \[
	\stochastic x 
	\bigl(\exists y/\{x\}\bigr)
	\stochastic z
	\biggl[\,
		(x = y = z) 
		\lor
		\Bigl[
		\neg(z = x \not = y) 
		\land 
		\neg(z = y \not= x) 
		\land 
		\bigl(\exists y/\{x\}\bigr)\,
		x = y\,
		\Bigr]
	\biggr],
    \]
and consider the semantic game $G\bigl(\{1, 2, 3\}, \varphi_{\mathrm{MH}\scriptsize{\stochastic}}'\bigr)$ in which Nature follows a behavioral strategy $\lambda$ that assigns values to $x$ and $z$ according to a uniform probability distribution. 

Although he no longer assigns values to $x$ and $z$, \abelard\ still plays a role in the semantic game because of the subformula 
    \[
        \neg(x \not= y = z) \land\neg(x  = z \not= y) \land \bigl(\exists y/\{x\}\bigr)\,x = y,
    \] 
which we treat as a ternary conjunction (denoted $\psi$ below). Let \(h_{abc} = \bigl(\varphi_{\mathrm{MH}\scriptsize{\stochastic}}', (x,a), (y,b), (z,c)\bigl)\), 
let  $\tau$ be the pure strategy for \abelard\ defined by 
    \begin{align*}
	\tau(h_{abc}\!\concat\psi) 
	&= 
	\begin{cases}
		\neg(z = x \not = y)		& \text{if \(c = a \not= b\)}, \\
		\neg(z = y \not= x)		& \text{if \(c = b \not= a\)}, \\
		\bigl(\exists y/\{x\}\bigr)\,x = y	& \text{otherwise},
	\end{cases}
    \end{align*}
and let $\nu$ be the mixed strategy for \abelard\ such that \(\nu(\tau) = 1\).

After the initial values of $x$, $y$, and $z$ have been set,  
\eloise\ should choose the left disjunct if and only if all three values are the same.
Let \(h_{a} = \bigl(\varphi_{\mathrm{MH}\scriptsize{\stochastic}}', (x,a)\bigr)\), 
and define \(\sigma_{b} \in S_{\exists}\) by 
\begin{align*}
	\sigma_{b}(h_{a}) &= (y,b), \\[6pt]
	\sigma_{b}(h_{abc}) &= 
		\begin{cases}
			\,x = y = z
			& 
			\text{if \(a = b = c\),} 
			\\
			\psi
			&
			\text{otherwise,}
		\end{cases} \\
	\sigma_{b}\Bigl((h_{abc}\!\concat\psi)\concat\bigl(\exists y/\{x\}\bigr)\,x = y\Bigr) &= 
		\min\bigl(\{1,2,3\}\setminus\{c\}\bigr).
\end{align*}
Define $\sigma^{b} \in S_{\exists}$ similarly except that 
    \(
	\sigma^{b}\Bigl((h_{abc}\!\concat\psi)\concat\bigl(\exists y/\{x\}\bigr)\,x = y\Bigr) = 
		\max\bigl(\{1,2,3\}\setminus\{c\}\bigr).
    \)
Let \(\mu \in \Delta(S_{\exists})\) be the mixed strategy according to which \(\mu(\sigma_{b}) = 1/6 = \mu(\sigma^{b})\).
(Strictly speaking, the formulas \(\neg(x \not= y = z)\) and \(\neg(x = z \not= y)\) abbreviate disjunctions, but we will treat them as literals.) 


\begin{figure}[htbp]
	\centering
		\resizebox{4in}{!}{
		\begin{tikzpicture}
			\node (root) 	{};
			\node (b)			[rounded corners,draw] 
										[right=of root, xshift=10mm]	{$\exists y/\{x\}$}
											edge[very thick] node[auto,swap] {1} (root);
			\node (ba)		[circle,draw] 
										[above right=of b, yshift=50mm] {$\stochastic z$}
											edge[very thick] node[auto,swap] {1} (b);
			\node (baa)		[circle,fill=gray]
										[above right=of ba, yshift=20mm, label=right:{$\exists$}]{}
											edge node[auto,swap] {1} (ba);
			\node (bab)	[rounded corners,draw]
										[right=of ba]{$\exists y/\{x\}$}
											edge[very thick] node[auto,swap] {2} (ba);
			\node (baba)	[circle,fill] 
										[above right=of bab, xshift=-0.5mm, yshift=-4mm, label=right:{$\exists$}]{}
											edge[very thick] node[auto,swap] {1} (bab);
			\node (babc)	[circle,fill] 
										[below right=of bab, xshift=-0.5mm, yshift=4mm, label=right:{$\forall$}]{}
											edge[very thick] node[auto] {3} (bab);
			\node (bac)	[rounded corners,draw]
										[below right=of ba, yshift=-20mm]{$\exists y/\{x\}$}
											edge[very thick] node[auto] {3} (ba);
			\node (baca)	[circle,fill] 
										[above right=of bac, xshift=-0.5mm, yshift=-4mm, label=right:{$\exists$}]{}
											edge[very thick] node[auto,swap] {1} (bac);
			\node (bacc)	[circle,fill] 
										[below right=of bac, xshift=-0.5mm, yshift=4mm, label=right:{$\forall$}]{}
											edge[very thick] node[auto] {2} (bac);
			\node (bb)		[circle,draw] 
										[right=of b, xshift=45mm]	{$\stochastic z$}
											edge[very thick] node[auto,swap] {2} (b);
			\node (bba)		[circle,fill=gray]
										[above right=of bb, yshift=20mm, label=right:{$\forall$}]{}
											edge node[auto,swap] {1} (bb);
			\node (bbb)	[circle,fill=gray]
										[right=of bb, label=right:{$\forall$}] {}
											edge node[auto,swap] {2} (bb);
			\node (bbc)	[rounded corners,draw]
										[below right=of bb, yshift=-20mm]{$\exists y/\{x\}$}
											edge[very thick] node[auto] {3} (bb);
			\node (bbca)	[circle,fill] 
										[above right=of bbc, xshift=-0.5mm, yshift=-4mm, label=right:{$\exists$}]{}
											edge[very thick] node[auto,swap] {1} (bbc);
			\node (bbcc)	[circle,fill] 
										[below right=of bbc, xshift=-0.5mm, yshift=4mm, label=right:{$\forall$}]{}
											edge[very thick] node[auto] {2} (bbc);
			\node (bc)		[circle,draw] 
										[below right=of b, yshift=-50mm]{$\stochastic z$}
											edge[very thick] node[auto] {3} (b);
			\node (bca)		[circle,fill=gray]
										[above right=of bc, yshift=20mm, label=right:{$\forall$}]{}
											edge node[auto,swap] {1} (bc);
			\node (bcb)	[rounded corners,draw]
										[right=of bc]{$\exists y/\{x\}$}
											edge[very thick] node[auto,swap] {2} (bc);
			\node (bcba)	[circle,fill] 
										[above right=of bcb, xshift=-0.5mm, yshift=-4mm, label=right:{$\exists$}]{}
											edge[very thick] node[auto,swap] {1} (bcb);
			\node (bcbc)	[circle,fill] 
										[below right=of bcb, xshift=-0.5mm, yshift=4mm, label=right:{$\forall$}]{}
											edge[very thick] node[auto] {3} (bcb);
			\node (bcc)	[circle,fill=gray]
										[below right=of bc, yshift=-20mm, label=right:{$\forall$}] {}
											edge node[auto] {3} (bc);
		\end{tikzpicture}
		}
			\caption{A portion of the semantic game  for \(\varphi_{\mathrm{MHS}}'\)
			}
  \label{Figure-PhiMHprime2GameTree}
\end{figure}

Figure \ref{Figure-PhiMHprime2GameTree} shows the histories that occur when \(x = 1\), \abelard\ follows $\nu$, and \eloise\ follows $\mu$. If initially \(y = 1\) and \(z = 1\), \eloise\ wins by choosing the left disjunct \(x = y = z\). 
If instead \(z = 2\) or \(z = 3\), 
then \eloise\ will choose the right disjunct $\psi$, 
after which \abelard\ will pick the right conjunct. 
\eloise\ will then reset the value of $y$ to be distinct from the value of $z$, 
winning when she sticks with her initial choice. 
If initially \(y = 2\) or \(y = 3\),
then \eloise\ will choose $\psi$, 
after which \abelard\ will win, if possible, by picking the appropriate conjunct based on the value of $z$. 
Otherwise, he will pick \(\bigl(\exists y/\{x\}\bigr)\,x = y\), 
in which case \eloise\ will reset the value of $y$ to be distinct from the value of $z$, 
but this time she will lose when she sticks with her original choice and win when she switches.

If we restrict our attention to those histories, 
highlighted in Figure \ref{Figure-PhiMHprime2GameTree}, 
in which 
Nature selects a value of $z$ that differs from the values of both $x$ and $y$, 
we can see that \eloise\ wins one-half of the time when she sticks with her original $y$-value,
and wins one-half of the time when she switches.

To conclude this section, let us consider the following stochastic variant $\varphi_{\mathrm{MH}\scriptsize{\stochastic}}$ of the original Monty Hall sentence,%
\footnote{
Sandu considers a variant of $\varphi_{\mathrm{MH}}$ in which only the first quantifier is stochastic \cite[page 247]{Sandu:2015}.}
\[
	\stochastic x 
	\bigl(\exists y/\{x\}\bigr)
	\stochastic z
	\Bigl[\,
		(z \not= x \land z \not= y) 
		\implies 
		\bigl(\exists y/\{x\}\bigr)\,
		x = y\,
	\Bigr], 
\]
which is an abbreviation for 
\[
	\stochastic x 
	\bigl(\exists y/\{x\}\bigr)
	\stochastic z
	\Bigl[\,
		(z = x \lor z = y)
		\lor 
		\bigl(\exists y/\{x\}\bigr)\,
		x = y\,
	\Bigr].
\]
Once again we assume that Nature follows a behavioral strategy $\lambda$ that assigns values to $x$ and $z$ according to a uniform probability distribution. Since the above sentence has no conjunctions or universal quantifiers, \abelard\ plays no role in its semantic game.
Let \(\sigma_{b}, \sigma_{b}' \in S_{\exists}\) be defined as before, with \eloise\ initially setting the value of $y$ to be $b$, then sticking or switching, respectively. Let \(\mu \in \Delta(S_{\exists})\) be the mixed strategy according to which \(\mu(\sigma_{b}) = 1/6 = \mu(\sigma_{b}')\). 

Unlike \abelard, Nature has no qualms about setting the value of $z$ equals to the value of $x$ or $y$. In fact, if Nature follows $\lambda$ and \eloise\ follows $\mu$, then \eloise\ will win by choosing the disjunct $(z = x \lor z = y)$ in 5/9 of all plays. In the remaining plays, exactly half will satisfy \(x = y\). Thus \eloise\ will win an additional 2/9 of all plays by correctly guessing the value of $x$, regardless of whether she sticks with her initial guess or switches. Thus, the truth value of $\varphi_{\mathrm{MH}\scriptsize\stochastic}$ is 7/9.


\section{Sleeping Beauty}
\label{Section:SleepingBeauty}
We now turn our attention to a related puzzle, called the Sleeping Beauty problem, popularized by Adam Elga:
%
\begin{quotation}
\noindent
Some researchers are going to put you to sleep.
During the two days that your sleep will last, 
they will briefly wake you up either once or twice,
depending on the toss of a fair coin 
(Heads: once; Tails: twice).
After each waking, 
they will put you back to sleep 
with a drug that makes you forget that waking.

When you are first awakened, to what degree ought you believe that the outcome of the coin toss is Heads?
\cite[page 143]{Elga:2000}
\end{quotation}
Two answers, 1/2 and 1/3, have been defended in the literature. 
The proponents of each answer are known as halfers and thirders, respectively.

Elga argues in favor of 1/3.
Suppose the first waking occurs on Monday, 
the second on Tuesday.
Then, immediately after waking, 
you will be in one of three possible situations:
$H_{1}$, 
the coin landed Heads and it is Monday;
$T_{1}$, 
the coin landed Tails and it is Monday; or
$T_{2}$,
the coin landed Tails and it is Tuesday.
According to Elga, all three are equally likely. 
For suppose that after waking you are told that the coin landed Tails. 
At that moment, you know that you are in situation $T_{1}$ or $T_{2}$, 
but\,---\,because of the drug\,---\,you cannot tell which.
Furthermore, 
you have no reason to suspect that $T_{1}$ is more or less likely than $T_{2}$, 
given that the coin landed tails.
Consequently, 
\[
	P(T_{1}\mid T_{1}\cup T_{2}) = P(T_{2}\mid T_{1}\cup T_{2}),
\]
which implies \(P(T_{1}) = P(T_{2})\)
\cite[page 144]{Elga:2000}.

Now suppose instead that, 
after waking you up, 
the researchers inform you that it is Monday.
Then you would know that you are in situation $H_{1}$ or $T_{1}$, 
which of the two being determined by the toss of a fair coin.
In fact, 
it doesn't matter whether the researchers toss the coin before or after waking you the first time. 
(You are woken on Monday regardless of the outcome.)
So you might suppose that the researchers have yet to toss the coin, 
in which case  
$P(H_{1}\mid H_{1}\cup T_1)$ 
is simply the probability that a fair coin that has yet to be tossed will land Heads.
Hence
\(P(H_{1}\mid H_{1}\cup T_{1}) = 1/2\). 
A simple calculation then shows that 
\(P(H_{1}) = P(T_{1})\).
Therefore 
\[
	P(H_{1}) = P(T_{1}) = P(T_{2}),
\] 
which implies \(P(H_{1}) = 1/3\)
\cite[page 145]{Elga:2000}.

David Lewis argues, 
contra Elga, 
that you do not gain any information relevant to the outcome of a fair coin toss simply by being awake 
since you knew all along that you would be awakened at least once during the experiment.
Lewis agrees with Elga that \(P(T_{1}) = P(T_{2})\), 
but challenges his assertion that 
\(P(H_{1}\mid H_{1}\cup T_{1}) = 1/2\). 
Instead, 
Lewis claims that one should assign equal prior probabilities to Heads and Tails:
\[
	P(H_{1}) = 1/2 = P(T_{1} \cup T_{2}).
\]
It follows that 
\(P(T_{1}) = 1/4 = P(T_{2})\).
Consequently,
when the researchers inform you that it is Monday, 
you should update your beliefs accordingly:
\[
P(H_{1}\mid H_{1}\cup T_{1}) = 
\frac{P(H_{1})}{P(H_{1}) + P(T_{1})} = 2/3.
\]
Lewis points out that he and Elga both agree that probability of Heads depends on what day it is, since 
\(P(H_{1}\mid H_{1}\cup T_{1}) = P(H_{1}) + 1/6\), 
but they disagree as to whether 
the prior or 
the posterior probability  
should be taken to be 1/2
\cite{Lewis:2001}.

Cian Dorr comes to Elga's defense 
by considering a variation of the experiment 
in which the researchers have two amnesia-inducing drugs instead of just one 
\cite{Dorr:2002}.
If the coin lands Tails,
they will administer the same drug as before, 
so that your experience is exactly the same as in the original experiment. 
If the coin lands Heads, however,
they will administer a weaker drug
so that, when you wake up the second time, 
your memories of the first awakening will be delayed for exactly one minute.
Thus, 
immediately after being awoken, 
it is possible that you are in the situation $H_{2}$,
the coin landed heads and it is Tuesday.
Since your subjective experience in each of the four situations is identical,
you should consider them all to be equally likely:
\[
	P(H_{1}) = P(H_{2}) = P(T_{1}) = P(T_{2}) = 1/4.
\]
After one minute passes, 
either you will remember being woken up on Monday, 
or you will not.
If you do, 
your memories will confirm that you are in the situation $H_{2}$. 
If you do not, 
your failure to remember (drug-induced or not) 
is evidence that you are currently experiencing  $H_{1}$, $T_{1}$, or $T_{2}$, 
at which point the probability that the coin landed Heads is
\[
	P(H_{1}\mid H_{1} \cup T_{1} \cup T_{2}) =
	\frac{P(H_{1})}{P(H_{1}) + P(T_{1}) + P(T_{2})} = 
	1/3.
\]
Arntzenius independently developed a similar variant of the Sleeping Beauty problem in which the subject of the experiment is a vivid dreamer who can distinguish wake from dream by pinching herself
\cite[pages 363--364]{Arntzenius:2003}.   
It is worth noting that whether Dorr's (and presumably Arntzenius's) variations are analogous to the original version of the problem is controversial
\cite{Bradley:2003,Horgan:2004}.

\bigskip\bigskip

We are now ready to formalize the Sleeping Beauty problem using stochastic IF logic.
At the moment she is awakened,
Beauty thinks to herself: 
``Given that I am awake, a fair coin must have been tossed, but I don't know whether it landed Heads or Tails. 
Furthermore, because of the amnesia-inducing drug I may have been given, I am unsure whether it is Monday or Tuesday.''
To help determine what her credence should be, she decides to model her predicament using the following stochastic IF sentence $\varphi_{\mathrm{SB}}$,
\[
	\stochastic x
	\stochastic t
			\Bigl[
				\Awake(x,t) \implies 
					\bigl(
						\Heads(x) \hor{\{x,t\}} \Tails(x)
					\bigr)
			\Bigr],
\]
interpreted in a structure $\str{M}$ with universe $\{1,2\}$ and equipped with the following relations: 
\begin{align*}
	\Heads^{\str{M}} &= \{1\} = \Monday^{\str{M}}\!\!\!\!\!\!, \\
	\Tails^{\str{M}} &= \{2\} = \Tuesday^{\str{M}}\!\!\!\!\!\!, \\
	\Awake^{\str{M}} &= \bigl\{(1,1), (2,1), (2,2)\bigr\}.
\end{align*}
Here $x$ represents the result of the coin toss, 
and $t$ represents the current time.

Next Beauty analyzes the semantic game 
$G(\str{M},\varphi_{\mathrm{SB}})$, 
which begins with Nature selecting values for $x$ and $t$.
Recall that the implication inside the square brackets is an abbreviation for
\[
	\neg\Awake(x,t) \lor \bigl(\Heads(x) \hor{\{x,t\}} \Tails(x)\bigr).
\]
Thus,
\eloise\ must first choose between
$\neg\Awake(x,t)$ 
and 
\(
	\bigl(\Heads(x) \hor{\{x,t\}} \Tails(x)\bigr).
\)
If she chooses the left disjunct, 
the game ends.
If she chooses the right disjunct, 
she must choose between the atomic formulas
$\Heads(x)$ or $\Tails(x)$
without knowing the value of $x$ or $t$.

Considering the players' possible strategies, 
Beauty reasons that Nature should follow a behavioral strategy $\lambda$ according to which 
\(\lambda(x \colonequals 1) = 1/2 = \lambda(x \colonequals 2)\) 
since the coin is fair. 
Beauty further assumes that 
\(\lambda(t \colonequals 1) = 1/2 = \lambda(t \colonequals 2)\)
since there is no reason for Nature to prefer one day over the other.

Assuming $\lambda$ is fixed, how should \eloise\ play?
After the partial history 
\(h_{ab} = \bigl(\varphi_{\mathrm{SB}}, (x,a), (t,b)\bigr)\),
she is confronted with the unslashed disjunction
\[
	\neg\Awake(x,t) \lor \bigl(\Heads(x) \hor{\{x,t\}} \Tails(x)\bigr).
\]
\eloise\ should clearly choose the left disjunct if and only if 
\((a,b) = (1,2)\), 
so let
\begin{align*}
	\sigma_{1}(h_{ab}) &=
		\begin{cases}
			\neg\Awake(x,t)		& 
				\text{if \((a,b) = (1,2)\)}, \\
			\bigl(\Heads(x) \hor{\{x,t\}} \Tails(x)\bigr) & \text{otherwise},
		\end{cases} \\
	\sigma_{1}\Bigl(h_{ab}\concat\bigl(\Heads(x) \hor{\{x,t\}} \Tails(x)\bigr)\Bigr) &=
		\Heads(x),
\end{align*}
and define $\sigma_{2}$ similarly except that
\(
	\sigma_{2}\Bigl(h_{ab}\concat\bigl(\Heads(x) \hor{\{x,t\}} \Tails(x)\bigr)\Bigr) = \Tails(x).
\)
All of her other strategies are dominated, 
so let $\mu_{p}$ be the mixed strategy over 
$\{\sigma_{1},\, \sigma_{2}\}$ 
such that 
\(\mu_{p}(\sigma_{1}) = p\)
and
\(\mu_{p}(\sigma_{2}) = 1 - p\).


\begin{figure}[htbp]
	\centering
		\begin{tikzpicture}
			\node (root) 
            			[circle,draw] {$\stochastic x$};
			\node (1)			
            			[circle,draw] 
            			[above right=of root, yshift=20mm] {$\stochastic t$}
            			edge node[auto,swap] {$\frac{1}{2}(\mathrm{Heads})$} (root);
			\node (11)			
            			[rounded corners,draw] 
            			[above right=of 1, xshift=5mm] {$\lor_{\!/\{x,t\}}$}
            			edge node[auto,swap] {$\frac{1}{2}(\mathrm{Mon})$} (1);
			\node (111)			
            			[circle,fill] 
            			[above right=of 11, xshift=5mm, label=right:{$\exists$}] {}
            			edge node[auto,swap] {$p$} (11);
			\node (112)			
            			[circle,fill] 
            			[below right=of 11, xshift=5mm, label=right:{$\forall$}] {}
            			edge node[auto] {$1 - p$} (11);
			\node (12)			
            			[circle,fill=gray] 
            			[below right=of 1, label=right:{$\exists$}] {}
            			edge node[auto,swap] {$\frac{1}{2}(\mathrm{Tue})$} (1);
			\node (2)			
            			[circle,draw] 
            			[below right=of root, yshift=-20mm] {$\stochastic t$}
            			edge node[auto] {$\frac{1}{2}(\mathrm{Tails})$} (root);
			\node (21)			
            			[rounded corners,draw] 
            			[above right=of 2, xshift=5mm, yshift=5mm] {$\lor_{\!/\{x,t\}}$}
            			edge node[auto,swap] {$\frac{1}{2}(\mathrm{Mon})$} (2)
            			edge [dotted] node[auto] {} (11);
			\node (211)			
            			[circle,fill] 
            			[above right=of 21, xshift=5mm, label=right:{$\forall$}] {}
            			edge node[auto,swap] {$p$} (21);
			\node (212)			
            			[circle,fill] 
            			[below right=of 21, xshift=5mm, label=right:{$\exists$}] {}
            			edge node[auto] {$1 - p$} (21);
			\node (22)			
            			[rounded corners,draw] 
            			[below right=of 2, xshift=5mm, yshift=-5mm] {$\lor_{\!/\{x,t\}}$}
            			edge node[auto] {$\frac{1}{2}(\mathrm{Tue})$} (2)
            			edge [dotted] node[auto] {} (21);
			\node (221)			
            			[circle,fill] 
            			[above right=of 22, xshift=5mm, label=right:{$\forall$}] {}
            			edge node[auto,swap] {$p$} (22);
			\node (222)			
            			[circle,fill] 
            			[below right=of 22, xshift=5mm, label=right:{$\exists$}] {}
            			edge node[auto] {$1 - p$} (22);
		\end{tikzpicture}
			\caption{The semantic game for \(\varphi_{\mathrm{SB}}\)
			}
  \label{Figure-PhiSBGameTree}
\end{figure}

As shown in Figure \ref{Figure-PhiSBGameTree}, 
when \((a,b) = (1,1)\),
\eloise\ 
wins if she follows $\sigma_{1}$ and
loses if she follows $\sigma_{2}$.
If \((a,b) = (1,2)\),
\eloise\ will win regardless of which strategy she follows.
If \(a = 2\), 
\eloise\ wins if she follows $\sigma_{2}$ and loses if she follows $\sigma_{1}$.
Hence
\eloise's expected utility from $\mu_{p}$ is
\begin{align*}
	U(\lambda,\mu_{p}) &= 
	\frac{p}{4} + 
	\frac{1}{4} +
	\frac{1 - p}{2} = 
	\frac{3 - p}{4}.
\end{align*}
Thus, \eloise\ wins exactly half of the plays in which she follows 
$\sigma_{1}$ 
and three-quarters of the plays in which she follows 
$\sigma_{2}$.
Hence her optimal strategy is to always follow $\sigma_{2}$.
Therefore \(\str{M} \models^{3/4}\!\varphi_{\mathrm{SB}}(\lambda)\).

However, 
Beauty is not so much interested in the truth value of the sentence $\varphi_{\mathrm{SB}}$ relative to $\lambda$ 
as she is in \eloise's chances of winning when  
\((a,b) \in \Awake^{\str{M}}\)\!\!\!\!\!.\;\,
If we exclude those histories in which \((a,b) = (1,2)\)
(grayed out in Figure \ref{Figure-PhiSBGameTree}),
then \eloise's chances of winning become
\[
	\frac
	{\frac{p}{4} + \frac{1 - p}{2}}
	{\frac{3}{4}} = 
	\frac{2 - p}{3}.
\]
Thus, 
\eloise\ wins one-third of the histories in which she chooses $\Heads(x)$ 
and two-thirds of the histories in which she chooses $\Tails(x)$.
Therefore, Beauty concludes that the probability that the coin landed Heads is 1/3.%
\footnote{
    Beauty's reasoning most closely matches Horgan's argument \cite{Horgan:2004} 
    that one should assign each of the four hypotheses 
    $H_{1}$, $H_{2}$, $T_{1}$, and $T_{2}$ a prior probability of 1/4. 
    See also Hitchcock's diachronic Dutch Book argument \cite{Hitchcock:2004}.
}


\begin{figure}[htbp]
	\centering
		\begin{tikzpicture}
			\node (root) 
            			[circle,draw] {$\stochastic x$};
			\node (1)			
            			[circle,draw] 
            			[above right=of root, yshift=20mm] {$\forall t$}
            			edge node[auto,swap] {$\frac{1}{2}(\mathrm{Heads})$} (root);
			\node (11)			
            			[rounded corners,draw] 
            			[above right=of 1, xshift=5mm] {$\lor_{\!/\{x,t\}}$}
            			edge node[auto,swap] {$1(\mathrm{Mon})$} (1);
			\node (111)			
            			[circle,fill] 
            			[above right=of 11, xshift=5mm, label=right:{$\exists$}] {}
            			edge node[auto,swap] {$p$} (11);
			\node (112)			
            			[circle,fill] 
            			[below right=of 11, xshift=5mm, label=right:{$\forall$}] {}
            			edge node[auto] {$1 - p$} (11);
			\node (12)			
            			[circle,fill=gray] 
            			[below right=of 1, label=right:{$\exists$}] {}
            			edge node[auto,swap] {$0(\mathrm{Tue})$} (1);
			\node (2)			
            			[circle,draw] 
            			[below right=of root, yshift=-20mm] {$\forall t$}
            			edge node[auto] {$\frac{1}{2}(\mathrm{Tails})$} (root);
			\node (21)			
            			[rounded corners,draw] 
            			[above right=of 2, xshift=5mm, yshift=5mm] {$\lor_{\!/\{x,t\}}$}
            			edge node[auto,swap] {$q(\mathrm{Mon})$} (2)
            			edge [dotted] node[auto] {} (11);
			\node (211)			
            			[circle,fill] 
            			[above right=of 21, xshift=5mm, label=right:{$\forall$}] {}
            			edge node[auto,swap] {$p$} (21);
			\node (212)			
            			[circle,fill] 
            			[below right=of 21, xshift=5mm, label=right:{$\exists$}] {}
            			edge node[auto] {$1 - p$} (21);
			\node (22)			
            			[rounded corners,draw] 
            			[below right=of 2, xshift=5mm, yshift=-5mm] {$\lor_{\!/\{x,t\}}$}
            			edge node[auto] {$(1 - q)(\mathrm{Tue})$} (2)
            			edge [dotted] node[auto] {} (21);
			\node (221)			
            			[circle,fill] 
            			[above right=of 22, xshift=5mm, label=right:{$\forall$}] {}
            			edge node[auto,swap] {$p$} (22);
			\node (222)			
            			[circle,fill] 
            			[below right=of 22, xshift=5mm, label=right:{$\exists$}] {}
            			edge node[auto] {$1 - p$} (22);
		\end{tikzpicture}
			\caption{The semantic game for \(\varphi_{\mathrm{SB}}'\)
			}
  \label{Figure-PhiSBprimeGameTree}
\end{figure}

A  halfer might object that Beauty made a mistake by considering the wrong stochastic IF sentence.
Instead of $\varphi_{\mathrm{SB}}$,
Beauty should have considered the following alternative sentence $\varphi_{\mathrm{SB}}'$,
\[
	\stochastic x
	\forall t
			\Bigl[
				\Awake(x,t) \implies 
					\bigl(
						\Heads(x) \hor{\{x,t\}} \Tails(x)
					\bigr)
			\Bigr],
\]
in which the universal player picks the value of $t$ instead of Nature. 
If she had, she would have found that 
Nature's behavioral strategy 
 \(\lambda(x \colonequals 1) = 1/2 = \lambda(x \colonequals 2)\) is stipulated by the fact that the coin is fair,
and that \eloise's pure and mixed strategies are the same as before.
Since \eloise\ wins any history in which \((a,b) = (1,2)\), 
\abelard\ should always set \(t \colonequals 1\) when \(x \colonequals 1\). 
Let $\nu_{q}$ be the mixed strategy for \abelard\ such that 
    \begin{align*}
	\nu_{q}(x \colonequals 1, t \colonequals 1) = 1 \qquad\text{and}\qquad&
	\nu_{q}(x \colonequals 1, t \colonequals 2) = 0, \\
	\nu_{q}(x \colonequals 2, t \colonequals 1) = q \qquad\text{and}\qquad&
	\nu_{q}(x \colonequals 2, t \colonequals 2) = 1- q.
    \end{align*}

\eloise's expected utility when 
she follows $\mu_{p}$, and
\abelard\ follows $\nu_{q}$ is
    \[
	U(\lambda, \mu_{p}, \nu_{q}) =
	\frac{p}{2} + 
	\frac{(1 - p)q}{2} + 
	\frac{(1 - p)(1 - q)}{2} = 
	1/2.
    \]
Notice that 
\(U(\lambda, \mu_{p}, \nu_{q})\) 
depends on neither $p$ nor $q$. 
Thus, \eloise\ wins exactly half of the time 
no matter how she guesses, reflecting the fact that the coin is fair.
Furthermore,
if we restrict our attention to those histories in which the value of $t$ is 1 (i.e., Monday), 
 \eloise's chances of winning become
    \[
	\frac{
            \frac{1}{2}p + 
            	\frac{1}{2}(1 - p)q 
	}{
		\frac{1}{2} + 
		\frac{1}{2}q
	} 
	=
	\frac{p + (1 - p)q}{1 + q}
    \]
When \(q = 0\), \eloise\ wins with probability $p$; 
when \(q = 1/2\), she wins with probability \((p + 1)/3\); 
and when \(q = 1\), she wins with probability 1/2.
This shows that, if the researchers always wake Beauty on Monday when the coin lands Heads, and always wake her on Tuesday (and only Tuesday) when the coin lands Tails, then, immediately after waking, her credence that the coin landed Heads should be 1/2, but, after being told that it is Monday, she becomes certain that the coin landed Heads.
At the other extreme, if the researchers always wake Beauty on Monday (and never on Tuesday), regardless of the result of the coin toss, then learning that it is Monday gives her no information about the coin toss, so her credence should remain 1/2.
Finally, suppose the researchers wake Beauty exactly once during the experiment: on Monday if the coin lands Heads, and on Monday \emph{or} Tuesday (with equal probability) if the coin lands Tails. Then her initial credence that the coin landed Heads should be 1/2, but should increase to 2/3 after she learns that it is Monday.

In summary, when \abelard\ follows $\nu_{1/2}$ and \eloise\ follows $\sigma_{1}$, she will win one-half of all plays, and
two-thirds of the plays in which the value of $t$ is 1. This is the correct result according to halfers. 
Unfortunately for them, 
\abelard's strategy $\nu_{1/2}$ \ is dominated by $\nu_{1}$, so it would be irrational for him to follow it. 


\begin{figure}[htbp]
	\centering
		\begin{tikzpicture}
			\node (root) 
            			[circle,draw] {$\stochastic x$};
			\node (1)			
            			[circle,draw] 
            			[above right=of root, yshift=20mm] {$\stochastic t$}
            			edge node[auto,swap] {$\frac{1}{2}(\mathrm{Heads})$} (root);
			\node (11)			
            			[rounded corners,draw] 
            			[above right=of 1, xshift=5mm] {$\lor_{\!/\{x,t\}}$}
            			edge node[auto,swap] {$1(\mathrm{Mon})$} (1);
			\node (111)			
            			[circle,fill] 
            			[above right=of 11, xshift=5mm, label=right:{$\exists$}] {}
            			edge node[auto,swap] {$p$} (11);
			\node (112)			
            			[circle,fill] 
            			[below right=of 11, xshift=5mm, label=right:{$\forall$}] {}
            			edge node[auto] {$1 - p$} (11);
			\node (12)			
            			[circle,fill=gray] 
            			[below right=of 1, label=right:{$\exists$}] {}
            			edge node[auto,swap] {$0(\mathrm{Tue})$} (1);
			\node (2)			
            			[circle,draw] 
            			[below right=of root, yshift=-20mm] {$\stochastic t$}
            			edge node[auto] {$\frac{1}{2}(\mathrm{Tails})$} (root);
			\node (21)			
            			[rounded corners,draw] 
            			[above right=of 2, xshift=5mm, yshift=5mm] {$\lor_{\!/\{x,t\}}$}
            			edge node[auto,swap] {$\frac{1}{2}(\mathrm{Mon})$} (2)
            			edge [dotted] node[auto] {} (11);
			\node (211)			
            			[circle,fill] 
            			[above right=of 21, xshift=5mm, label=right:{$\forall$}] {}
            			edge node[auto,swap] {$p$} (21);
			\node (212)			
            			[circle,fill] 
            			[below right=of 21, xshift=5mm, label=right:{$\exists$}] {}
            			edge node[auto] {$1 - p$} (21);
			\node (22)			
            			[rounded corners,draw] 
            			[below right=of 2, xshift=5mm, yshift=-5mm] {$\lor_{\!/\{x,t\}}$}
            			edge node[auto] {$\frac{1}{2}(\mathrm{Tue})$} (2)
            			edge [dotted] node[auto] {} (21);
			\node (221)			
            			[circle,fill] 
            			[above right=of 22, xshift=5mm, label=right:{$\forall$}] {}
            			edge node[auto,swap] {$p$} (22);
			\node (222)			
            			[circle,fill] 
            			[below right=of 22, xshift=5mm, label=right:{$\exists$}] {}
            			edge node[auto] {$1 - p$} (22);
		\end{tikzpicture}
			\caption{The semantic game for \(\varphi_{\mathrm{SB}}\) when Nature follows $\lambda'$
			}
  \label{Figure-PhiSBprime2GameTree}
\end{figure}
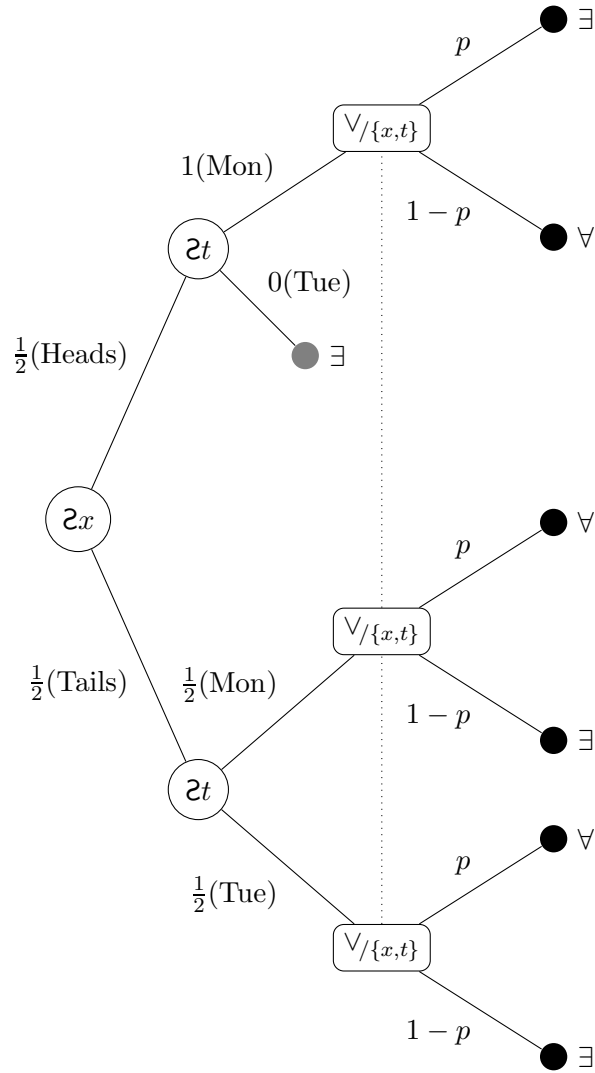

The preceding analysis suggests another possibility, however. Instead of proposing the alternative stochastic IF sentence $\varphi_{\mathrm{SB}}'$, a halfer might assert that, in the semantic game for $\varphi_{\mathrm{SB}}$, Nature should follow the alternative behavioral strategy $\lambda'$ according to which 
    \begin{align*}
	\lambda'(x \colonequals 1, t \colonequals 1) = 1/2 \qquad\text{and}\qquad&
	\lambda'(x \colonequals 1, t \colonequals 2) = 0, \\
	\lambda'(x \colonequals 2, t \colonequals 1) = 1/4 \qquad\text{and}\qquad&
	\lambda'(x \colonequals 2, t \colonequals 2) = 1/4.
    \end{align*}
However, Beauty's prior analysis of $\varphi_{\mathrm{SB}}'$ shows that the game depicted in Figure \ref{Figure-PhiSBprime2GameTree} does not accurately formalize her predicament. Rather, it formalizes an experiment in which the subject is awakened exactly once: on Monday if the coin lands Heads, and on Monday \emph{or} Tuesday (with equal probability) if the coin lands Tails. Thus, we are forced to conclude that Lewis' argument in favor of 1/2 is correct for this variant of the Sleeping Beauty problem, 
but is incorrect when applied to the original version.


\section{Conclusion}

The dual goals of the present article have been to deepen our understanding of the Monty Hall and Sleeping Beauty problems by formalizing them in a logical language, while simultaneously using these famous problems to introduce a natural extension of first-order logic with imperfect information.

We began by showing how the original Monty Hall problem and a variant in which the host is not required to offer the contestant the opportunity to switch doors could be viewed as semantic games of IF sentences, where \abelard\ plays the role of the host and \eloise\ plays the role of the contestant. We then modeled two further variants of the Monty Hall problem in which the host is indifferent to the outcome as semantic games with chance moves. 

We also showed how Beauty could use semantic games with chance moves to analyze her predicament, leading to the conclusion that the thirders are correct. Finally, we explained how semantic games with chance moves could be used to demonstrate that Lewis's argument in favor of 1/2 applies to a variant of the Sleeping Beauty problem. Although we doubt  ours will be the last word on the matter, we hope that future contributors to the Sleeping Beauty debate will follow our example by presenting explicitly their formalizations of the problem.

\bibliographystyle{abbrv}
\bibliography{main}

\end{document}